\documentclass[journal]{IEEEtran}

\ifCLASSINFOpdf
\else
\fi
\usepackage{amsmath}
\usepackage{amssymb}
\usepackage[dvipdfmx]{graphicx}
\usepackage{color}
\usepackage{dsfont}
\usepackage{cite}

\usepackage{theorem}
\theorembodyfont{\rmfamily}

\newtheorem{definition}{Definition}

\newtheorem{theorem}{Theorem}
\newtheorem{lemma}{Lemma}
\newtheorem{proposition}{Proposition}

\renewcommand{\IEEEQED}{{\hfill \IEEEQEDopen}\\}

\newcommand{\bm}[1]{\boldsymbol{#1}}

\DeclareMathOperator{\sfsin}{\mathsf{sin}}
\DeclareMathOperator{\sfcos}{\mathsf{cos}}

\DeclareMathOperator{\sfdiag}{\mathsf{diag}}

\DeclareMathOperator{\sfdet}{\mathsf{det}}

\DeclareMathOperator{\sfRe}{\mathsf{Re}}
\DeclareMathOperator{\sfIm}{\mathsf{Im}}

\DeclareMathOperator{\sfker}{\mathsf{ker}}

\DeclareMathOperator{\sfspan}{\mathsf{span}}

\DeclareMathOperator{\sfint}{\mathsf{int}}

\newcommand{\mat}[1]{\left[\: \begin{matrix} #1 \end{matrix} \:\right]}
\newcommand{\spliteq}[1]{\begin{split} #1 \end{split}}
\newcommand{\simode}[1]{\left\{\:  \begin{aligned} #1 \end{aligned} \right.}

\begin{document}
\title{
Equilibrium-Independent Passivity of Power Systems: A Link Between Classical and\\
Two-Axis Synchronous Generator Models
}

\author{
        Takayuki~Ishizaki,~\IEEEmembership{Member,~IEEE,}
        Taku~Nishino,
        and~Aranya Chakrabortty,~\IEEEmembership{Senior Member,~IEEE}% <-this % stops a space
\thanks{
T.~Ishizaki and T.~Nishino are with 
%T.~Ishizaki is with 
		Graduate School of Engineering, Tokyo Institute of Technology, Meguro,
		Tokyo, Japan. 
		e-mail: \{ishizaki, ta.nishino\}@lim.sc.e.titech.ac.jp.}% <-this % stops a space
\thanks{A. Chakrabortty is with 
		the Electrical and Computer Engineering Department, North Carolina State University, 
		Raleigh, NC 27695 USA.
		e-mail: achakra2@ncsu.ed.}% <-this % stops a space
\thanks{This research was supported by JST-Mirai Program 18077648.}
\thanks{A preliminary version of this paper is found as \cite{ishizaki2021necessity}.}
\thanks{Manuscript received xxxx; revised xxxx.}}

\markboth{Journal of \LaTeX\ Class Files,~Vol.~X, No.~X, XXXX}%
{Shell \MakeLowercase{\textit{et al.}}: Bare Demo of IEEEtran.cls for IEEE Journals}

\maketitle

\begin{abstract}
We study the equilibrium-independent (EI) passivity of a nonlinear power system composed of two-axis generator models. 
The model of our interest consists of a feedback interconnection of linear and nonlinear subsystems, called mechanical and electromagnetic subsystems. 
We mathematically prove the following three facts by analyzing the nonlinear electromagnetic subsystem. 
First, a lossless transmission network is necessary for the EI passivity of the electromagnetic subsystem. 
Second, the convexity of a strain energy function characterizes the largest set of equilibria over which the electromagnetic subsystem is EI passive.
Finally, we prove that the strain energy function for the network of the two-axis generator models is convex if and only if its flux linkage dynamics is stable, and the strain energy function for the network of the classical generator models derived by singular perturbation approximation of the flux linkage dynamics is convex.
Numerical simulation of the IEEE  9-bus power system model demonstrates the practical implications of the various mathematical results.
In particular, we validate that the convex domain of the strain energy function over which the electromagnetic subsystem is EI passive is almost identical to the set of all stable equilibria.
This result is also generalized to lossy power systems based on our finding that the convexity of the strain energy function is equivalent to the positive semidefiniteness of a synchronizing torque coefficient matrix.
\end{abstract}


\IEEEpeerreviewmaketitle

\section{Introduction}\label{sec:intro}

\IEEEPARstart Mathematical results from systems and control theory have made significant contributions to the modeling, stability analysis, and control of electric power systems \cite{kundur1994power}. 
In particular, passivity theory has been used extensively in the literature to understand various fundamental physical properties of power systems as well as the methods for controlling them by analyzing the dynamical models of synchronous electric generators. 
For example, \cite{arcak2007passivity} and \cite{yang2019distributed} have derived results on modeling and stability of power systems using the notion of passivity and differential passivity. 
Papers such as \cite{trip2016internal,stegink2017unifying,wang2018distributed,li2016connecting}, on the other hand, have presented various centralized and distributed control algorithms for automatic generation control (AGC) of synchronous generators, based on the passivity of the input-output map from their turbine mechanical power input to their frequency deviations. 
Alternative approaches to control design have been presented in papers such as \cite{ortega2005transient}, where the field voltage of the electrical excitation systems of the generators is used to passify feedback connections in the generator models using generator state measurements composed of different types of components.

The common thread in all of these papers is the analysis of the dynamical models of synchronous machines, which are the most commonly used electrical machines for transmission-level power supply. 
However, the complete dynamical model of a synchronous generator can be quite complex, as it consists of a large list of state variables arising from both stator and rotor flux linkages, currents and voltages, evolving over multiple time scales. 
For example, following Sections~5.1 and 5.2 in \cite{sauer2017power}, one can see that a typical multi-timescale model of a synchronous generator, excluding a voltage regulator and a governor, consists of at least nine states, including the states of the rotor angle, stator transient, field-winding flux linkage, and damper-winding flux linkage. 

The more than nine-dimensional model is usually reduced to a six-dimensional model as the dynamics of the stator transients are much faster than the other states, and the value of the network resistance is negligible compared to the network inductance at the transmission level.  
Furthermore, if two flux states of the damper-winding dynamics are ignored due to their relatively fast dynamics, the dimension of the model is further reduced by two states, resulting in a four-dimensional model, consisting of the two-dimensional electromechanical swing equations of the generator rotor, and the two-dimensional electromagnetic dynamics of the field winding and damper winding flux linkages. 
The resulting four-dimensional generator model is referred to as the \textit{two-axis model}, which serves as a good working model of a synchronous generator. 
It should be noted that most control designs take an additional step in simplifying the dynamics by ignoring the flux linkage of the damper winding, resulting in the so-called one-axis or flux-decay model \cite{sadamoto2019dynamic}. 
Simplified models such as the \textit{classical model}, which will be discussed in this paper, are also commonly used. 
See Chapter~5 in \cite{sauer2017power} for more details on the derivations of these different  models with different resolutions of dynamics.

In this paper, we consider the two-axis model as the starting point of our discussion, since it is the most generic dynamical model for practical applications of a synchronous machine at the transmission level. 
Using this model, we present mathematical analysis of the equilibrium-independent (EI) passivity of power systems, which allows to systematically determine the set of asymptotically stable equilibria, i.e., stable power flow distributions, in terms of kinetic and strain energies. 
The EI passivity \cite{hines2011equilibrium,simpson2019equilibrium}  is a relatively new notion of passivity defined as the passivity with respect to a set of feasible equilibria.
The entire power system model considered in this paper is represented as the negative feedback of linear and nonlinear subsystems, referred to as the \textit{mechanical subsystem} and the \textit{electromagnetic subsystem}, respectively.
In particular, because proving the EI passivity of the mechanical subsystem is straightforward as it is the parallel connection of first-order stable linear systems, we focus only on proving the EI passivity of the nonlinear electromagnetic subsystem. 
We use the Bregman divergence of a strain energy function as the storage function for this electromagnetic subsystem. 
Note that the existing results on passivity analysis of power systems such as in \cite{trip2016internal, stegink2017unifying, wang2018distributed, yang2019distributed,DePersis2017} only consider the case of  one-axis or classical models of the generators.

We further elaborate the necessary conditions for the EI passivity of the electromagnetic subsystem, and its relation to the EI passivity of the classical model network. 
In particular, we mathematically prove that
\begin{itemize}
\item[(a)] a lossless transmission network is not only sufficient, but also necessary for the electromagnetic subsystem to be EI passive,
\item[(b)] the convex domain of the strain energy function is equal to the ``largest" set of equilibria over which the electromagnetic subsystem is EI passive, and 
%\item[(c)] the EI passivity of the two-axis model is equivalent to the stability of its \red{flux linkage dynamics}, and the EI passivity of a classical model derived by singular perturbation approximation of the flux linkage dynamics.
\item[(c)] the strain energy function of the two-axis model network is convex if and only if its flux linkage dynamics is stable, and the strain energy function of a classical model network derived by singular perturbation approximation (SPA) of the flux linkage dynamics is convex.
\end{itemize}
An additional finding from these analyses is that the convexity of the strain energy function is equivalent to the positive semidefiniteness of a \textit{synchronizing torque coefficient matrix} that is generalized to multi-machine power systems.
Note that the conventional concept of synchronizing torque coefficients is well defined only for the single-machine infinite-bus model with the classical generator model \cite{sauer2017power}.

%These novel findings are proven using linearization and Kron reduction (a procedure to transform the DAE representation of the power system model to an equivalent ODE representation \cite{dorfler2013kron, ishizaki2018graph}) of the generators models.
%In their proofs, we utilize several remarkable properties of an admittance matrix derived from the original DAE representation before the Kron reduction.

Several remarks are in order. 
Most of the existing works on passivity analysis of power systems assume that the transmission network is lossless. 
Based on this premise, the nonlinear dynamics from the mechanical power input to the angular frequency deviation of the generators are shown to be passive.
As exceptions, \cite{ortega2005transient} and \cite{yang2019distributed} have addressed the case of lossy power systems. 
The former proposes a feedback passification method for power systems with lossy transmission lines, while the latter proposes an approximate storage function based on the notion of numerical energy functions \cite{chang1995direct}.
A related work \cite{narasimhamurthi1984existence} shows that a classical energy function cannot be used as a Lyapunov function for lossy power systems.
However, the ``necessity" of such losslessness is not easy to prove by nonlinear analysis because a storage function for passivity is not uniquely determined, implying that some other storage functions may exist to prove passivity. 
Our necessity analysis also contributes to the analysis of the largest set of equilibria and the discovery of a new link between the two-axis and classical generator models in terms of the EI passivity.

On the other hand, it should be noted that even if losslessness is essential for the EI passivity, it may not be necessary  for the stability of power systems, or more specifically, for the small-signal stability of power flow distributions. 
This is because the set of equilibria over which a system is EI passive is in general only a subset of the stable equilibria.
Nevertheless, it is demonstrated in this paper that
\begin{itemize}
\item[(d)]  the largest set of equilibria over which the electromagnetic subsystem is EI passive, or equivalently the convex domain of its strain energy function, is almost identical to the set of all stable equilibria.
\end{itemize}
We validate this fact through numerical simulations of the IEEE 9-bus power system model \cite{sauer2017power}, considering both lossless and lossy transmission.
For the stability analysis of lossy power systems, the convexity of the strain energy function should be replaced with the positive semidefiniteness of the aforementioned synchronizing torque coefficient matrix.
This practical finding indicates the rationality of discussing the stability of even lossy power systems in terms of the EI passivity concept.

The remainder of this paper is organized as follows.
In Section~\ref{sec:sysdes}, we describe a power system model composed of the two-axis generator models.
We also explain how the classical model can be derived as a special case of the two-axis model.
In Section~\ref{sec:eipass}, we prove that the nonlinear electromagnetic subsystem is EI passive if the transmission network is lossless.
In Section~\ref{sec:anakronlin}, based on linearization, we elaborate on the necessity of the lossless transmission for the EI passivity.
Section~\ref{sec:conv} finds a new link between the two-axis and classical generator models in terms of the EI passivity. 
Section~\ref{sec:numex} presents numerical simulations. 
Finally, some concluding remarks are made in Section~\ref{sec:conc}.

\vspace{2mm}
\noindent \textbf{Notation}~ 
%The notation used in this paper is generally standard.
We denote the set of real values by $\mathbb{R}$, 
the set of non-negative real values by $\mathbb{R}_{\geq 0}$, 
the sphere by $\mathbb{S}$,
the diagonal matrix whose diagonal entries are $\{a_1,\ldots,a_N\}$ by $\sfdiag(a_i)_{i\in \{1,\ldots,N\}}$ or simply by $\sfdiag(a_i)$,
the imaginary unit by $\bm{j}$,
the real and imaginary parts of a complex number $\bm{y}$ by $\sfRe[\bm{y}]$ and $\sfIm[\bm{y}]$, respectively,
the all-ones vector by $\mathds{1}$,
the subspace spanned by a vector $v$ by $\sfspan \{v\}$,
the interior of a set $\mathcal{D}$ by $\sfint \mathcal{D}$,
the positive definiteness and semidefiniteness of a symmetric matrix $A$ by $A\succ 0$ and $A\succeq 0$, respectively,
the negative definiteness and semidefiniteness by the converse symbols, 
%the kernel of a matrix $A$ by $\sfker A$,
the element-wise complex conjugate of a complex matrix $Z $ by $\overline{Z}$,
the Jacobian of a vector field $f:\mathbb{R}^n \rightarrow \mathbb{R}^m$ evaluated at $x^{\star}$ by $\tfrac{\partial f}{\partial x}(x^{\star})$,
the gradient of a function $U:\mathbb{R}^n \rightarrow \mathbb{R}$ by
\[
\nabla U(x) :=
\mat{
\tfrac{\partial U}{\partial x_1}(x) & \cdots & \tfrac{\partial U}{\partial x_n}(x)
},
\]
the Bregman divergence of $U$ by
\[
\mathfrak{B}_{x^{\star} } \bigl[
U(x) 
\bigr]:=
U(x) - U(x^{\star}) - \nabla U (x^{\star}) (x-x^{\star}),
\]
%\[
%\mathfrak{B}[U^{\rm red}](x;x^{\star}) :=
%U(x) - U(x^{\star}) - \nabla U (x^{\star}) (x-x^{\star}),
%\]
and the Hessian of $U$ by %$\nabla^2U(x)$.
\[
\nabla^2 U(x):=
\mat{
\tfrac{\partial^2 U}{\partial x_1^2} (x) & \cdots & \tfrac{\partial^2 U}{\partial x_1 \partial x_n} (x) \\
\vdots & \ddots & \vdots \\
\tfrac{\partial^2 U}{\partial x_n \partial x_1} (x) & \cdots &\tfrac{\partial^2 U}{\partial x_n^2} (x)
}.
\]
A square matrix $L$ is said to be a weighted graph Laplacian if it is symmetric, its off-diagonal entries are non-positive, and every row sum is zero.
The main theoretical results are stated as Theorems, while not main but remarkable facts are stated as Propositions.
Lemmas are provided for technical purposes.

%%%%%%%%%%%%%%%%%%%%%%%%%%%%%%%%%%%%%%%%%%%%%%%%%%%%%%%%%%%%%%%%%%%%%%%%%%%%%%%%%%%%%%%%%%%%%%%%%%%%

\section{System Description}\label{sec:sysdes}
\subsection{Power System Model for Mathematical Analysis}\label{sec:DAEmod}

We review a power system model found in the literature \cite{sauer2017power,Anderson2008,sadamoto2019dynamic}.
For mathematical analysis, we consider a power system where one generator is connected to each of all buses.
In such a model consisting only of generators, some generators can be considered as induction motors that consume active power \cite{ishizaki2018graph}. 
In particular, these induction motors can also be considered as a load model with power-frequency droop control \cite{simpson2013synchronization, guggilam2017optimizing} that consumes specified active power when time constants are sufficiently small.
We will demonstrate it by numerical simulation in Section~\ref{sec:numex}.

%The procedure reported in Appendix~\ref{sec:redload} yields the following model, where constant impedance loads are involved as a part of transmission lines.
%We remark that constant power loads can also be involved in a similar way. We will demonstrate it by numerical simulation in Section~\ref{sec:numex}.

%In such a power system model consisting only of generators, some generators may consume active power depending on the setting of transmission line constants. 
%Those synchronous generators with negative power can be regarded as a type of loads equipped with grid-forming inverters, such as virtual synchronous generators \cite{lin2020research}.

\begin{figure}[t]
\centering
\includegraphics[width = .60\linewidth]{ 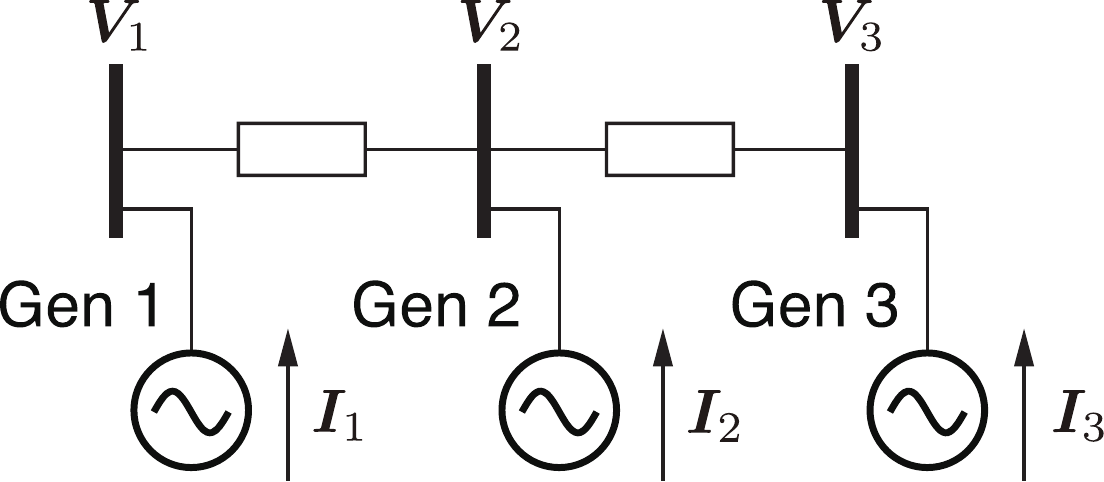}
%\medskip
\caption{Example of 3-bus power system.}
\label{fig3busgen}
\end{figure}

\smallskip
\subsubsection{Transmission Network Model}\label{sec:DAEmodel}

Consider a power system composed of $N$ buses.
See Fig.~\ref{fig3busgen} for an example of three buses.
Let $\bm{y}_{ij}$ and $c_{ij}$ denote the admittance and ground capacitance of the transmission line between Bus~$i$ and Bus~$j$. 
Then, the admittance matrix $\bm{Y}$ of the transmission network, each of whose line is represented by a $\pi$-type equivalent circuit, is given as
\begin{equation}
\bm{Y}_{ij}
=
\left\{
\begin{array}{cl}
-\bm{y}_{ij}, & i\neq j \\
\sum\limits_{j=1}^{N} \left(\bm{y}_{ij} + \bm{j} \frac{\omega_0 c_{ij}}{2}\right) , & i=j
\end{array}
\right.
\end{equation}
where $\bm{Y}_{ij}$ denotes the $(i,j)$-element of $\bm{Y}$.
The network equation of all buses is given as 
\begin{align}\label{eq:ohmY}
\mat{
  \bm{I}_1\\
  \vdots\\
  \bm{I}_N
}
 =
\mat{
  \bm{Y}_{11} & \cdots & \bm{Y}_{1N}\\
  \vdots & \ddots & \vdots\\
  \bm{Y}_{N1} & \cdots & \bm{Y}_{NN}
}
\mat{
  \bm{V}_1\\
  \vdots\\
  \bm{V}_N
}.
\end{align}
The real part and imaginary part of $\bm{Y}$, denoted as 
\begin{equation}\label{eq:GandB}
\bm{Y} = G + \bm{j} B,
\end{equation}
are referred to as the \textit{conductance matrix} and \textit{susceptance matrix}, respectively.

\smallskip
\subsubsection{Two-Axis Generator Model}

We consider a synchronous generator at Bus $i$.
Let $E_{{\rm d}i}\in \mathbb{R}$ denote the field-winding flux linkage, 
$E_{{\rm d}i}\in \mathbb{R}$ denote the damper-winding flux linkage, and
$\delta_i\in \mathbb{R}$ denote the rotor angle relative to the frame rotating at the system angular frequency $\omega_0$.
The bus voltage phasor of Bus $i$, denoted by $\bm{V}_i \in \mathbb{C}$, and the bus current phasor from the generator to Bus $i$, denoted by $\bm{I}_i \in \mathbb{C}$, are related by
\begin{subequations}\label{eq:genmodel}
\begin{equation}\label{eq:phVI}
\simode{
\bm{V}_{{\rm q}i} &= E_{{\rm q}i} - X_{i}' \bm{I}_{{\rm d}i} \\
\bm{V}_{{\rm d}i} &= E_{{\rm d}i} + X_{i}' \bm{I}_{{\rm q}i}
}
\end{equation}
where $X'_{i}$ is the transient reactance, and
\[
\spliteq{
\bm{V}_{{\rm d}i} &:=
|\bm{V}_i| \sfsin (\delta_i - \angle \bm{V}_i ) ,
 \\
\bm{V}_{{\rm q}i} &:=
|\bm{V}_i| \sfcos (\delta_i - \angle \bm{V}_i ) ,
}
\quad
\spliteq{
\bm{I}_{{\rm d}i} &:=
|\bm{I}_i| \sfsin (\delta_i - \angle \bm{I}_i ) ,
 \\
\bm{I}_{{\rm q}i} &:=
|\bm{I}_i| \sfcos (\delta_i - \angle \bm{I}_i ).
}
\]
By definition, $|\bm{V}_i| \in \mathbb{R}_{\geq 0}$, $|\bm{I}_i| \in \mathbb{R}_{\geq 0}$, $\angle \bm{V}_i \in \mathbb{S}$, and $\angle \bm{I}_i \in \mathbb{S}$.
The equation \eqref{eq:phVI} can be seen as an ``output equation" of the synchronous generator where $E_{{\rm q}i}$, $E_{{\rm d}i}$ and $\delta_i$ are the state variables, and $\bm{V}_i$ and $\bm{I}_i$ can be regarded as the interaction input and output signals, respectively.
Note that the active power and the reactive power, defined as
\[
P_i := \sfRe \left[ \bm{V}_i \overline{ \bm{I} }_i \right]
,\quad
Q_i := \sfIm \left[ \bm{V}_i \overline{ \bm{I} }_i \right],
\]
can be chosen as an equivalent interaction output, equal to
\begin{equation}\label{eq:PiQis}
P_i  =  \bm{V}_{{\rm q}i} \bm{I}_{{\rm q}i} +
\bm{V}_{{\rm d}i} \bm{I}_{{\rm d}i} ,\quad
Q_i  =  
\bm{V}_{{\rm q}i} \bm{I}_{{\rm d}i} -
\bm{V}_{{\rm d}i} \bm{I}_{{\rm q}i}.
\end{equation}

Next, we consider the dynamics of the synchronous generator, composed of a two-dimensional swing equation and a two-dimensional flux linkage dynamics.
Let $ \omega_i \in \mathbb{R}$ be the angular frequency deviation relative to the standard value $\omega_0$, and $P_{{\rm m}i}\in \mathbb{R}$ be the mechanical input.
Then, the synchronous generator dynamics is given as
\begin{equation}\label{eq:gendynVIst}
\simode{
\dot{\delta}_i&= \omega_0   \omega_i\\
M_i    \dot{\omega}_i&= 
 - D_i \omega_i  
 - P_i 
+P_{{\rm m}i} 
\\
\tau_{{\rm d}i} \dot{E}_{{\rm q}i} & = 
- \tfrac{ X_{i} }{ X_{i}' }E_{{\rm q}i}
+\left(
\tfrac{ X_{i} }{ X_{i}' } -1
\right)
\bm{V}_{{\rm q}i}
+ V_{{\rm fd}i}^{\star} \\
\tau_{{\rm q}i} \dot{E}_{{\rm d}i} & = 
- \tfrac{ X_{i} }{ X_{i}' }E_{{\rm d}i}
+\left(
\tfrac{ X_{i} }{ X_{i}' } -1
\right)
\bm{V}_{{\rm d}i}
}
\end{equation}
\end{subequations}
where, 
$M_i$ is the inertia constant,
$D_i$ is the damping coefficient, 
$\tau_{{\rm d}i} $ and $\tau_{{\rm q}i} $ are the time constants of the flux linkage dynamics,
$X_{i} $ is the synchronous reactance, and 
$V_{{\rm fd}i}^{\star}$ is the field voltage, supposed to be a positive constant throughout this paper because it is much faster than the mechanical input.
All constants here are positive.
%Note that $P_{{\rm m}i}$ and $V_{{\rm fd}i}$ are considered as slow and fast external inputs for generator control, respectively.
For standard parameters \cite{sauer2017power,sadamoto2019dynamic}, the transient reactance $X_{i}'$ is smaller than or equal to the synchronous reactance $X_{i}$.

Eliminating $\bm{I}_{{\rm d}i}$ and $\bm{I}_{{\rm q}i}$ from \eqref{eq:phVI} and \eqref{eq:PiQis}, we see that
\begin{equation}\label{eq:PiQis2}
\spliteq{
P_i & =  \frac{E_{{\rm q}i}}{X_{i}' } \bm{V}_{{\rm d}i}
- \frac{E_{{\rm d}i}}{X_{i}' } \bm{V}_{{\rm q}i}
+
\left(
\frac{1}{X_{i}' } - \frac{1}{X_{i}'} 
\right)
\bm{V}_{{\rm d}i} \bm{V}_{{\rm q}i}
, \\
Q_i & =  
\frac{E_{{\rm q}i}}{X_{i}' } \bm{V}_{{\rm q}i}
+\frac{E_{{\rm d}i}}{X_{i}' } \bm{V}_{{\rm d}i}
-
\left(
\frac{\bm{V}_{{\rm d}i}^2}{X_{i}' } 
+
\frac{\bm{V}_{{\rm q}i}^2 }{X_{i}'} 
\right),
}
\end{equation}
where  $\bm{V}_i$ can be regarded as the interaction input signal.
We remark that the resultant power system is represented as a differential algebraic equation (DAE) model where the dynamical generators in \eqref{eq:genmodel} are interconnected by the algebraic equation in \eqref{eq:ohmY}.
Furthermore, in a standard parameter setting, the ground capacitances are small enough for
\begin{equation}\label{eq:assumgc}
\underbrace{\textstyle
\sum\limits_{j=1}^N
\tfrac{\omega_0c_{ij}}{2}
}_{\beta_i}
 X_i' <1
 ,\quad \forall i \in \{1,\ldots,N\}
\end{equation}
to hold with the transient reactances.
We assume this condition because it is sufficient to guarantee the rationality of the algebraic constraint.
See Appendix~\ref{sec:supmath} for details.

%Recalling the standard definition in power systems engineering, we introduce the following terminology 
%\cite{kundur1994power, sauer2017power}.
%
%
%\begin{definition}
%The two-axis synchronous generator model in \eqref{eq:genmodel} is said to be a \textit{non-salient pole type} if 
%\begin{equation}\label{eq:nonsal}
%X_{{\rm d}i} = X_{{\rm q}i}
%,\quad
%X_{{\rm d}i}' = X_{{\rm q}i}',
%\end{equation}
%whose values are denoted by $X_i$ and $X_i'$, respectively.
%\end{definition}
%
%
%Similarly, the generator model is referred to as a salient pole type if \eqref{eq:nonsal} is not satisfied.
%Salient pole types, or projecting pole types, are mostly used in hydro power plants, while non-salient pole types, or round rotor types, are used in nuclear, gas, and thermal power plants.

\subsection{Approximation to Classical Generator Model}

In this subsection, we review the existing fact that the classical generator model can be derived as a special case of the two-axis generator model \cite{sauer2017power}.
Let us consider the case where the time constants $\tau_{{\rm q}i}$ and $\tau_{{\rm d}i}$ in the flux linkage dynamics are sufficiently small.
In particular, applying the SPA to the dynamics of $E_{{\rm q}i}$ and $E_{{\rm d}i}$ yields
\[
\simode{
0 & = 
- \tfrac{ X_{i} }{ X_{i}' }E_{{\rm q}i}
+\left(
\tfrac{ X_{i} }{ X_{i}' } -1
\right)
\bm{V}_{{\rm q}i}
+ V_{{\rm fd}i}^{\star} \\
0 & = 
- \tfrac{ X_{i} }{ X_{i}' }E_{{\rm d}i}
+\left(
\tfrac{ X_{i} }{ X_{i}' } -1
\right)
\bm{V}_{{\rm d}i}.
}
\]
Substituting this into \eqref{eq:phVI}, we have
\begin{subequations}\label{eq:clmodel}
\begin{equation}\label{eq:phVIcl}
\simode{
\bm{V}_{{\rm q}i} &= V_{{\rm fd}i}^{\star} - X_{i} \bm{I}_{{\rm d}i} \\
\bm{V}_{{\rm d}i} &= X_{i} \bm{I}_{{\rm q}i}.
}
\end{equation}
Then, the generator dynamics in \eqref{eq:gendynVIst} is reduced to the simple swing equation
\begin{equation}\label{eq:gendynVIstcl}
\simode{
\dot{\delta}_i&= \omega_0   \omega_i\\
M_i    \dot{\omega}_i&= 
 - D_i \omega_i  
 - P_i 
+P_{{\rm m}i} ,
}
\end{equation}
\end{subequations}
and the output equation in \eqref{eq:PiQis} can be simplified as
\[
P_i = 
\frac{V_{{\rm fd}i}^{\star}  }{X_i } \bm{V}_{{\rm d}i} 
,\quad
Q_i = 
\frac{V_{{\rm fd}i}^{\star}  }{X_i }  \bm{V}_{{\rm q}i}
- \frac{|\bm{V}_i|^2 }{X_i }.
\]

For the following discussion, we summarize the relation between the two-axis and classical generator models in the following proposition.

\begin{proposition}\label{prop:2axcl}
The two-axis generator model in \eqref{eq:genmodel} coincides with the classical generator model in \eqref{eq:clmodel} in the limit where both $\tau_{{\rm q}i}$ and $\tau_{{\rm d}i}$ are sufficiently small.
\end{proposition}

Note that the one-axis generator model \cite{sadamoto2019dynamic} can also be derived by a similar SPA.
Because of this relationship between the generator models, their steady state distributions are shown to be equivalent to each other.
In particular, the steady state distributions of the bus voltage phasors, the bus current phasors, and the rotor angles are identical for all power system models composed of the two-axis, one-axis, and classical generator models.
This is also true when different generator models are mixed.
Such a mixed system can be understood as a special case of the two-axis model network with different parameters.
Note that the transient behavior and stability properties of these different generator models are generally different because the time constants of the flux linkage dynamics are not very small.
One of the main results of this paper is to find a new link between the two-axis and classical generator models in terms of the largest set of stable equilibria based on the analysis of EI passivity, on the top of this relationship over the different timescale.

\section{EI Passivity of Lossless Power Systems}\label{sec:eipass}
\subsection{Review of EI Passivity}

Consider an input-affine nonlinear system
\begin{equation}\label{eq:nlsig}
\Sigma: \simode{
\dot{x} &= f(x) + Bu + Rd^{\star}\\
y &= h(x) + Du
}
\end{equation}
where $x \in \mathbb{R}^{n} $ is the state, $u \in \mathbb{R}^m$ is the input, and $y \in \mathbb{R}^m$ is the output.
Furthermore, $B \in \mathbb{R}^{n\times m}$, 
$D \in \mathbb{R}^{m\times m}$, and 
$R \in \mathbb{R}^{n\times p}$ 
are matrices,  
$d^{\star}\in \mathbb{R}^p$ is a constant input, and 
$f:\mathbb{R}^{n} \rightarrow \mathbb{R}^{n}$ and 
$h:\mathbb{R}^{n} \rightarrow \mathbb{R}^{m}$ are assumed to be sufficiently smooth, and satisfy 
\[
f(0)=0
,\quad
h(0)=0.
\]
Its state space, input space, and output space are denoted by $\mathcal{X}$, $\mathcal{U}$, and $\mathcal{Y}$, respectively.
Then, we introduce the notion of EI passivity \cite{hines2011equilibrium,simpson2019equilibrium}
or shifted passivity \cite{monshizadeh2019conditions} as follows.

\begin{definition}\label{def:eipassive}
Consider a nonlinear system $\Sigma$ in \eqref{eq:nlsig}.
Denote the set of feasible equilibria by
\begin{equation}\label{eq:asbleq}
\hspace{-22mm}
\mathcal{E} := 
\left\{
x^{\star}  \in \mathcal{X}: 
0 = f(x^{\star})+B u^{\star} + Rd^{\star}, \,
\exists u^{\star}  \in \mathcal{U}
\right\}.
\hspace{-15mm}
\end{equation}
Then, $\Sigma$ is said to be \textit{EI passive} over the equilibrium set $\mathcal{E}$ if for every $x^{\star} \in \mathcal{E}$, there exists a differentiable positive definite storage function $W_{x^{\star}}:\mathcal{X} \rightarrow \mathbb{R}_{\geq 0}$ such that $W_{x^{\star}} (x^{\star})=0$ and 
\begin{align}\label{eq:eiconpv}
\frac{d}{dt} W_{x^{\star}} \bigl( x(t) \bigr) \leq (u-u^{\star})^{\sf T} (y-y^{\star})
\end{align}
for all $u \in \mathcal{U}$ and $t \geq 0$, where $u^{\star} \in \mathcal{U}$ and $y^{\star} \in \mathcal{Y}$ are the constant input and outputs corresponding to $x^{\star} \in \mathcal{X}$.
In particular, $\Sigma$ is said to be \textit{strictly EI passive} if there additionally exists  a constant $\rho >0$ such that
\begin{align}\label{eq:eiconosp}
\frac{d}{dt} W_{x^{\star}} \bigl( x(t) \bigr) \leq (u-u^{\star})^{\sf T} (y-y^{\star})
-\rho\|y-y^{\star}\|^2
\end{align}
for all $u \in \mathcal{U}$ and $t \geq 0$.
\end{definition}

%In what follows, we abbreviate ``equilibrium-independent" as ``EI."
It has been shown in \cite{simpson2019equilibrium} that for every EI passive system, a storage function in the form of
\begin{equation}\label{eq:Brest}
W_{x^{\star}}(x) := 
\mathfrak{B}_{x^{\star} } \bigl[
U(x) 
\bigr]
%W_{x^{\star}}(x) = W(x) - W(x^{\star}) - \nabla W^{\sf T} (x^{\star}) (x-x^{\star})
\end{equation}
can be found.
This particular structure clarifies that the positive definiteness of the storage function $W_{x^{\star}}$ is relevant to the convexity of the energy function $U$.
Note that the EI passivity is reduced to the standard notion of the passivity if a feasible steady state configuration $(x^{\star},u^{\star},y^{\star})$ is specified.

\subsection{Equivalent ODE Representation of Power System}\label{sec:feedrep}

In the following, we derive an equivalent ordinary differential equation (ODE) model from the DAE model in Section~\ref{sec:sysdes}.
See Appendix~\ref{sec:kron} for the mathematical details of this derivation.
To this end, we define a \textit{reduced admittance matrix} as 
\begin{equation}\label{eq:Yred}
\bm{Y}^{\rm red}:= 
-\bm{j} 
\Bigl\{
%\underbrace{
\sfdiag \left( X_i' \right) 
-  \bm{j} \sfdiag \left( X_i' \right) \overline{\bm{Y}} \sfdiag \left( X_i' \right)
%}_{\bm{\varGamma}}
\Bigr\}^{-1},
\end{equation}
whose real and imaginary parts are denoted by
\begin{equation}\label{eq:GBred}
\bm{Y}^{\rm red} = G^{\rm red} + \bm{j} B^{\rm red}.
\end{equation}
In the following, $G^{\rm red}$ and $B^{\rm red}$ are referred to as the \textit{reduced conductance matrix} and \textit{reduced susceptance matrix}, respectively.
Furthermore, for
\[
\delta_{ij}:=\delta_i - \delta_j,
\]
we define the trigonometric functions
\begin{equation}\label{eq:defkh}
\spliteq{
k_{ij}(\delta_{ij})
&:=
-B_{ij}^{\rm red}
\sfcos\delta_{ij}
-
G_{ij}^{\rm red}
\sfsin\delta_{ij},
\\
h_{ij}(\delta_{ij})
&:=
-B_{ij}^{\rm red}
\sfsin\delta_{ij}
+
G_{ij}^{\rm red}
\sfcos\delta_{ij}
}
\end{equation}
where $G_{ij}^{\rm red}$ and $B_{ij}^{\rm red}$ denote the $(i,j)$-elements of $G^{\rm red}$ and $B^{\rm red}$, respectively.
Then, we obtain an equivalent ODE representation of the power system model, the  $i$th subsystem of which is given as
\begin{equation}\label{eq:kron2ax}
\simode{
\dot{\delta}_i&= \omega_0   \omega_i\\
M_i    \dot{\omega}_i&= 
- D_i \omega_i  
- P_i (z)
+ P_{{\rm m}i} 
\\
\tau_{{\rm d}i} \dot{E}_{{\rm q}i} & = 
- \tfrac{ X_i }{ X_i' } E_{{\rm q}i}
+\left(
X_i -X_i'
\right)
g_{{\rm q}i} (z)
+ V_{{\rm fd}i}^{\star} \\
\tau_{{\rm q}i} \dot{E}_{{\rm d}i} & = 
- \tfrac{ X_i }{ X_i' }E_{{\rm d}i}
+\left(
X_i -X_i'
\right)
g_{{\rm d}i} (z)
}
\end{equation}
where the active power output is obtained as
\begin{equation}\label{eq:deffi}
\spliteq{
P_i (z) :=
\sum_{j=1}^N
\Bigl[
&
\bigl\{ E_{{\rm q}i} h_{ij}(\delta_{ij}) - E_{{\rm d}i} k_{ij}(\delta_{ij}) \bigr\} E_{{\rm q}j} \\
+ & \bigl\{ E_{{\rm q}i} k_{ij}(\delta_{ij}) + E_{{\rm d}i} h_{ij}(\delta_{ij}) \bigr\} E_{{\rm d}j}
\Bigr],
}
\end{equation}
and the bus voltage phasor terms are obtained as
\begin{equation}\label{eq:kronterms}
\spliteq{
g_{{\rm q}i} (z) &:= \sum_{j=1}^N
\Bigl\{
k_{ij}(\delta_{ij}) E_{{\rm q}j}
- 
h_{ij}(\delta_{ij}) E_{{\rm d}j}
\Bigr\}, \\
g_{{\rm d}i} (z) &:= \sum_{j=1}^N
\Bigl\{
h_{ij}(\delta_{ij}) E_{{\rm q}j}
+
k_{ij}(\delta_{ij}) E_{{\rm d}j}
\Bigr\}.
}
\end{equation}
Note that $z$ denotes the vector composed of $\delta$, $E_{\rm q}$, and $E_{\rm d}$.

We consider the ODE equivalent as a feedback system of two subsystems.
One is the linear subsystem given as
\begin{equation}\label{eq:sys1}
\mathds{F}_i:
\simode{
M_i  \dot{\omega}_i&= 
- 
D_i
\omega_i  
 + 
v_i
\\
w_i &= \omega_0 \omega_i  ,
}
\end{equation}
which represents the mechanical dynamics of generators.
We denote the collection of $\mathds{F}_i$ by $\mathds{F}$.
In the following, this $\mathds{F}$ is referred to as the \textit{mechanical subsystem}.

\begin{figure}[t]
\centering
\includegraphics[width = .70\linewidth]{ 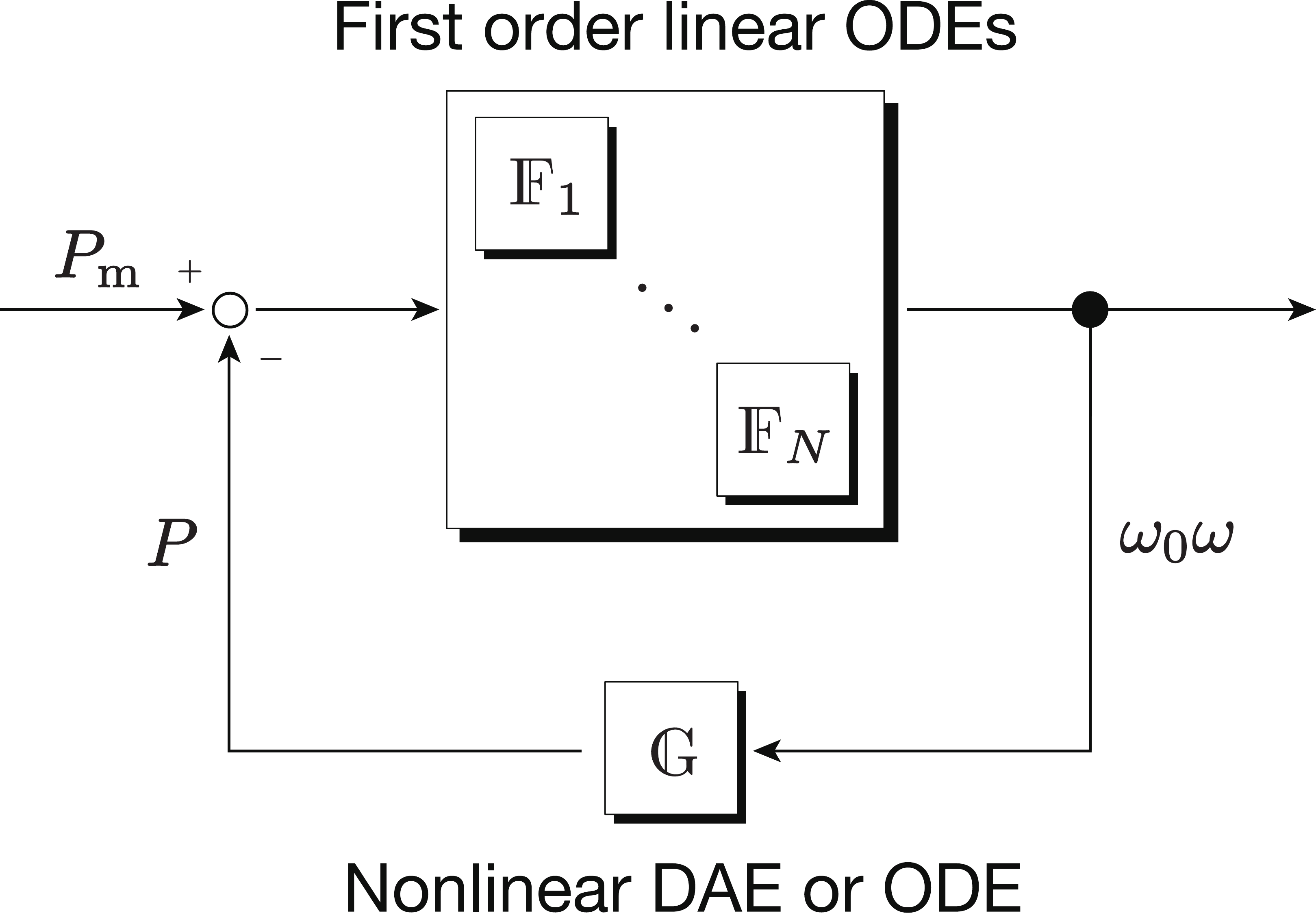}
%\medskip
\caption{Feedback representation of power system.}
\label{figblockFG}
\end{figure}

The other is the nonlinear subsystem given as
\begin{equation}\label{eq:sys2all}
\mathds{G}_i : 
\simode{
\dot{\delta}_i &= u_i 
\\
\tau_{{\rm d}i} \dot{E}_{{\rm q}i} & = 
- \tfrac{ X_i }{ X_i' } E_{{\rm q}i}
+\left(
X_i -X_i'
\right)
g_{{\rm q}i} (z)
+ V_{{\rm fd}i}^{\star} \\
\tau_{{\rm q}i} \dot{E}_{{\rm d}i} & = 
- \tfrac{ X_i }{ X_i' }E_{{\rm d}i}
+\left(
X_i -X_i'
\right)
g_{{\rm d}i} (z)
\\
y_i &= P_i(z).
}
\end{equation}
We denote the collection of \eqref{eq:sys2all} by $\mathds{G}$.
With a slight abuse of notation, we refer to it as the \textit{electromagnetic subsystem}, consisting of the flux linkage dynamics dependent on the bus voltage phasors in the transmission network.
The subsystems $\mathds{F}$ and $\mathds{G}$ are interconnected by the negative feedback given as
\begin{align}
v_i = P_{{\rm m}i}- y_i
,\quad
u_i = w_i
,\quad
\forall i \in \{1,\ldots,N\}.
\end{align}
This system representation is depicted in Fig.~\ref{figblockFG}, where the stacked versions of symbols are denoted by those without the subscript $i$, e.g., the stacked version of $P_{i}$ by $P$.

For the following discussion, we introduce the notion of losslessness as follows.

\begin{definition}
The transmission network is said to be \textit{lossless} if the reduced conductance matrix $G^{\rm red}$ in \eqref{eq:GBred} is zero.
\end{definition}

We remark that the reduced conductance matrix $G^{\rm red}$ is zero if and only if the conductance matrix $G$ is zero, meaning that there is no resistance in all transmission lines.
The sufficiency of this fact is relatively easy to prove, while the necessity is not so easy.
See Proposition~\ref{prop:GredG} in Appendix~\ref{sec:supmath} for the proof.

\subsection{EI Passivity of Electromagnetic Subsystem}

As a preliminary result of this paper, we prove that the power system composed of the two-axis generator models has the EI passivity if the transmission network is lossless.
To the best of the authors' knowledge, this finding on the EI passivity considering the two-axis generator model has not been reported in the literature.
This is a generalization of the result for the one-axis model, found as in \cite{wang2018distributed} for example.

Because $\mathds{F}$ is the parallel connection of the first-order stable linear systems, it is clear that $\mathds{F}$ is strictly EI passive if every $M_i$ and $D_i$ are positive.
On the other hand, the EI passivity of $\mathds{G}$ is not trivial.
We remark that, if $\mathds{G}$ is EI passive, then the entire power system model is strictly EI passive from the mechanical input $P_{\rm m}$ to the angular frequency deviation $\omega$.
In fact, when the transmission network is lossless, the EI passivity $\mathds{G}$ is proven as follows.

\begin{proposition}\label{prop:eqpassive}
Consider the electromagnetic subsystem $\mathds{G}$ in \eqref{eq:sys2all}.
Assume that the transmission network is lossless.
Define a strain energy function by
\begin{equation}\label{eq:defstEkron}
\spliteq{
U(z)   := &
 \sum_{i=1}^N 
\Biggl\{ 
\frac{X_i E_{{\rm q}i}^2 }{ 2 X_i '(X_i - X_i')}
+
\frac{X_i E_{{\rm d}i}^2 }{ 2 X_i '(X_i - X_i')} \\
  + &
\sum_{j=1}^N B^{\rm red}_{ij} \Bigl[
E_{{\rm q}i} \bigl(E_{{\rm q}j}  \sfcos \delta_{ij} 
-  E_{{\rm d}j} \sfsin \delta_{ij} \bigr) \\
& \hspace{10pt} + 
E_{{\rm d}i} \bigl(E_{{\rm d}j}  \sfcos \delta_{ij} 
+  E_{{\rm q}j}  \sfsin \delta_{ij} \bigr)
\Bigr]
\Biggr\},
}
\end{equation}
and the storage function by
\begin{equation}\label{eq:stfunW}
W_{z^{\star}}(z ) :=
\mathfrak{B}_{z^{\star}} \bigl[
U (z ) 
\bigr].
\end{equation}
Then, $\mathds{G}$ is EI passive over the equilibrium set 
\begin{equation}\label{eq:eqEG}
\mathcal{E} :=
{\sf int} \! \left\{
z^{\star}
\in \mathcal{Z}^{\star} :
\nabla^2 U(z^{\star}) \succeq 0
\right\}
\end{equation}
where $\mathcal{Z}^{\star} $ denotes the set of all feasible equilibria.
\end{proposition}

\begin{IEEEproof}
As long as $\mathcal{E}$ in \eqref{eq:eqEG} is not empty, there exists an open neighborhood of $z^{\star}$ over which $W_{z^{\star}}$ is positive definite for any feasible $z^{\star} \in \mathcal{E}$.
Note that the non-negativity of the storage function $W_{z^*}$ in \eqref{eq:stfunW} is written as
\[
U(z) 
\geq  U (z^{\star}) 
+ \nabla U (z^{\star})(z-z^{\star}).
\]
This is satisfied if both $z$ and $z^{\star}$ belong to the domain such that $U$ is convex.
It is known that the convexity of the twice-differentiable function $U$ is characterized by the positive semidefiniteness of its Hessian \cite{rockafellar1970convex}.

%Note that $\mathcal{E}$ denotes the set of all admissible steady states over which the twice-differentiable function $U$ is convex.
Using the state equation of $\mathds{G}$, we see that
the partial derivatives of $U $ are obtained as
\[
\spliteq{
\frac{\partial U }{\partial \delta_i} (z ) = P_i(z), \quad
\frac{\partial U }{\partial E_{{\rm q}i} }(z ) &= 
\frac{ V_{{\rm fd}i}^{\star}  -\tau_{{\rm d}i} \dot{E}_{{\rm q}i} }{ X_{i}-X_{i}'} \\
\frac{\partial U }{\partial E_{{\rm d}i} } (z) &= 
-\frac{ \tau_{{\rm q}i} \dot{E}_{{\rm d}i} }{ X_{i}-X_{i}'}.
}
\]
Therefore, along the trajectory of $\mathds{G}$, we have
\[
\spliteq{
\frac{d}{dt}
W_{z^{\star}} \bigl(z (t) \bigr)
& = \bigl\{
\nabla U  (z )
-
\nabla U  (z^{\star} )
\bigr\}
\dot{z} \\
& \leq (y - y^{\star})^{\sf T} (u - u^{\star}),
}
\]
where $u^{\star}$ must be zero for all feasible equilibria.
This proves that $\mathds{G}$ is EI passive over $\mathcal{E}$.
\end{IEEEproof}

Proposition~\ref{prop:eqpassive} shows that the electromagnetic subsystem $\mathds{G}$ is EI passive over the equilibrium set  $\mathcal{E}$ in \eqref{eq:eqEG}.
It should be emphasized that we cannot check how conservative such an equilibrium set $\mathcal{E}$ is.
This is because the storage function is not uniquely determined in general, and some other storage functions may exist to prove the EI passivity with respect to a larger equilibrium set.
%Furthermore, even the existence of $\mathcal{E}$ is not clear because the domain over which the storage function is positive definite is implicitly defined.
We will discuss these points in detail to show that it is indeed the ``largest."

\section{Necessity of Lossless Transmission for EI Passivity of Electromagnetic Subsystem}\label{sec:anakronlin}

In this section, we will prove that the transmission network must be lossless for the electromagnetic subsystem to be EI passive, i.e., Fact (a) in Section~\ref{sec:intro}.
To this end, we will elaborate on a linearized version, because the electromagnetic subsystem is EI passive over a set of equilibria only if its linearization is passive at each of all such equilibria.

%Throughout Sections~\ref{sec:anakronlin} and \ref{sec:conv}, we assume that every synchronous generator model is a non-salient pole type to make the Kron reduction, a procedure to transform the original DAE to an equivalent ODE, is mathematically tractable.
%The case of salient pole types will be studied numerically in Section~\ref{sec:numex}.
%The main objective here is to show that 
%\begin{itemize}
%\item the transmission network must be lossless for the electromagnetic subsystem to be EI passive, and
%\item the convex domain of a strain energy function captures the largest set of equilibria with respect to which the electromagnetic subsystem is EI passive
%\end{itemize}
%by analysis via the Kron reduction and linearization.

\subsection{Linearization}

In the rest of this paper, we denote the vector composed of $E_{\rm q}$ and $E_{\rm d}$ by $E$.
We consider linearizing the electromagnetic subsystem $\mathds{G}$ in \eqref{eq:sys2all}.
In particular, its linearized version is represented in the form 
\begin{equation}\label{eq:Gss}
{\mathds G}_{\vartriangle}: \simode{
\dot{{\delta}}_{\vartriangle}
&= \textstyle 
u_{\vartriangle}  \\
\tau  \dot{{E}}_{\vartriangle}
&= \textstyle 
A  E_{\vartriangle} + B \delta_{\vartriangle} \\
y_{\vartriangle} &= CE_{\vartriangle} + L \delta_{\vartriangle},
}
\end{equation}
%\end{subequations}
where $\delta_{\vartriangle}$ and $E_{\vartriangle}$ correspond to the deviations from the steady state values of $\delta^{\star}$ and $E^{\star}$, respectively, and
\[
\tau := \mat{
\sfdiag(\tau_{{\rm d}i}) & 0 \\
0 & \sfdiag(\tau_{{\rm q}i})
}.
\]
The system matrices $L \in \mathbb{R}^{N\times N}$, $A\in \mathbb{R}^{N\times N}$, $B\in \mathbb{R}^{2N\times N}$, and $C\in \mathbb{R}^{N\times 2N}$ will be specified in the following.

For the stacked version
\[
g(z):=\mat{g_{\rm q}(z) \\g_{\rm d}(z) }
\]
where $g_{\rm q}$ and $g_{\rm d}$ are, respectively, the stacked versions of $g_{{\rm q}i}$ and $g_{{\rm d}i}$ in \eqref{eq:kronterms}, we define 
\begin{subequations}\label{eq:sysmats}
\begin{equation}
\spliteq{
\hat{A}& := \frac{\partial g}{\partial E}(z^{\star})
+
I_2 \otimes \sfdiag \left(  \frac{X_i}{X_i'(X_i -X_i')}  \right), \\
\hat{B} & := \frac{\partial g}{\partial \delta}(z^{\star})
}
\end{equation}
where  $I_2$ denotes the two-dimensional identity matrix, and $\otimes$ denotes the Kronecker product.
Then, $A$ and $B$ are defined as
\begin{equation}
\spliteq{
A &:= \left\{ I_2 \otimes \sfdiag \left(  X_i -X_i'  \right) \right\} \hat{A},
\\
B &:= \left\{ I_2 \otimes \sfdiag \left(  X_i -X_i'  \right) \right\} \hat{B}.
}
\end{equation}
Furthermore, $L$ and $C$ are defined as
%\begin{equation}
%L:=\tfrac{\partial f}{\partial \delta}(z^{\star})
%,\quad
%C := \mat{
%\frac{\partial f}{\partial E_{\rm q}}(z^{\star}) &
%\frac{\partial f}{\partial E_{\rm d}}(z^{\star})
%}
%\end{equation}
\begin{equation}
L:=\frac{\partial P}{\partial \delta}(z^{\star})
,\quad
C := \frac{\partial P}{\partial E}(z^{\star})
\end{equation}
\end{subequations}
where $P$ is the stacked version of $P_i$ in \eqref{eq:deffi}.

\subsection{Necessity of Lossless Transmission for EI Passivity}

In this subsection, we prove that the transmission network must be lossless for the linearized electromagnetic subsystem ${\mathds G}_{\vartriangle}$ in \eqref{eq:Gss} to be passive, or equivalently, for its transfer matrix defined as
\begin{align}\label{eq:trGs}
H(s) :=  - \frac{1}{s} 
\underbrace{
\left\{ -C \bigl( s\tau  -A \bigr)^{-1} B - L \right\}
}_{\hat{H}(s)}
\end{align}
to be positive real.
For completeness, we present the standard definition of positive realness as follows.

\begin{definition}\label{def:trpf}
For a square transfer matrix $Q$, define
\begin{align}\label{eq:defOm0}
\Omega_0 := \left\{
\omega_0 \in \mathbb{R}: 
\mbox{ $\bm{j} \omega_0$ is a pole of $Q$}
\right\}.
\end{align}
Then, $Q$ is said to be \textit{positive real} if 
\begin{itemize}
\item all poles of $Q$ have non-positive real parts, \vskip 2pt
\item $Q(\bm{j} \omega) + Q^{\sf T}(-\bm{j} \omega) \succeq 0$ for all $\omega \in [0,\infty)\setminus \Omega_0$, and \vskip 2pt
\item every pure imaginary pole of $Q$ is simple, and its residue satisfies 
\begin{align*}
\lim_{s \rightarrow \bm{j} \omega_0} (s-\bm{j} \omega_0) Q(s) = 
\lim_{s\rightarrow \bm{j} \omega_0} \bigl\{ \overline{(s-\bm{j} \omega_0) \bm{j} Q(s)} \bigr\}^{\sf T} \! \succeq 0
\end{align*}
for all $\omega_0 \in \Omega_0$.
\end{itemize}
\end{definition}

To analyze the positive realness of $G$, the following dual notion will also be used \cite{petersen2010feedback}.

\begin{definition}\label{def:trni}
For a square transfer matrix $Q$ having no poles on the origin, define $\Omega_0$ in \eqref{eq:defOm0}.
Then, $Q$ is said to be \textit{negative imaginary} if 
\begin{itemize}
\item all poles of $Q$ have non-positive real parts, \vskip 2pt
\item $\bm{j}\{Q(\bm{j} \omega) - Q^{\sf T}(-\bm{j} \omega) \} \! \succeq  0$ for all $\omega \in \! (0,\infty)\setminus \Omega_0$, and \vskip 2pt
\item every pure imaginary pole of $Q$ is simple, and its residue satisfies 
\begin{align*}
\lim_{s \rightarrow \bm{j} \omega_0} (s-\bm{j} \omega_0) \bm{j} Q(s) = \!\!
\lim_{s\rightarrow \bm{j} \omega_0} \bigl\{ \overline{(s-\bm{j} \omega_0) \bm{j} Q(s)} \bigr\}^{\sf T} \! \succeq 0
\end{align*}
for all $\omega_0 \in \Omega_0$.
\end{itemize}
\end{definition}

In the rest of this paper, the matrix
\begin{equation}\label{eq:Lapmat}
L_0:= L-CA^{-1}B
\end{equation}
plays an important role.
The following theorem is one of the main mathematical findings in this paper.

\begin{theorem}\label{thm:EdynNI}
Consider the transfer matrices $H$ and $\hat{H}$ in \eqref{eq:trGs}.
For any feasible equilibrium $z^{\star}$ such that $\hat{H}$ is stable, $\hat{H}$ is negative imaginary if and only if the transmission network is lossless.
In addition, $H$ is positive real if and only if the transmission network is lossless, and
\begin{align}\label{eq:pdsp}
L_0 = L_0^{\sf T} \succeq 0,
\end{align}
where $L_0$ is defined as in \eqref{eq:Lapmat}.
\end{theorem}

\begin{IEEEproof}
First, we prove that if $\hat{H}$ is negative imaginary, then the transmission network is lossless, i.e., $G^{\rm red}$ is zero.
For $\hat{H}$ to be negative imaginary, $L$ must be symmetric as shown in \cite{petersen2010feedback}.
Calculating the difference between the $(i,j)$-element and $(j,i)$-element of $L$, we have 
\[
\spliteq{
\Delta L_{ij} := &
\tfrac{\partial P_i}{\partial \delta_j}(z^{\star}) -
\tfrac{\partial P_j}{\partial \delta_i}(z^{\star}) \\
= &
\bigl\{
h_{ij}(z^{\star}) + h_{ji}(z^{\star}) 
\bigr\}
(E_{{\rm d}i}^{\star} E_{{\rm q}j}^{\star} + E_{{\rm q}i}^{\star} E_{{\rm d}j}^{\star} ) 
 \\
- &
\bigl\{
k_{ij}(z^{\star}) - k_{ji}(z^{\star}) 
\bigr\}
(E_{{\rm q}i}^{\star} E_{{\rm q}j}^{\star} + E_{{\rm d}i}^{\star} E_{{\rm d}j}^{\star} )
}
\]
where $k_{ij}$ and $h_{ij}$ are defined as in \eqref{eq:defkh}.
Note that
\[
\spliteq{
h_{ij}(z^{\star}) + h_{ji}(z^{\star}) 
&= 2 G^{\rm red}_{ij} \sfcos \delta_{ij}^{\star} , \\
k_{ij}(z^{\star}) - k_{ji}(z^{\star}) 
&=- 2 G^{\rm red}_{ij} \sfsin \delta_{ij}^{\star}  .
}
\]
Thus, if $\Delta L_{ij}$ is zero for any admissible $z^{\star}$, then
\[
G^{\rm red}_{ij} = 0,\quad \forall i \neq j.
\]
In such a case, Proposition~\ref{prop:GredG} in Appendix~\ref{sec:supmath} proves that
\[
\sfdiag
\bigl\{
G^{\rm red}_{ii} (1 - \beta_i  X_i')
\bigr\}
\mathds{1}=0.
\]
On the premise of \eqref{eq:assumgc}, this equality means that every $G^{\rm red}_{ii}$ is also zero because $G^{\rm red}$ is positive semidefinite, or more specifically all diagonal elements of $G^{\rm red}$ are nonnegative, as shown in Proposition~\ref{prop:defBGred}.

Next, we prove that, if $G^{\rm red}$ is zero, then $\hat{H}$ is negative imaginary.
From the negative imaginary lemma \cite{petersen2010feedback}, we see that $\hat{H}$ is negative imaginary if $L$ is symmetric, and there exists a positive definite $V$ such that
\begin{align}\label{eq:cndQ}
\tilde{A}^{\sf T}V + V\tilde{A} \preceq 0
,\quad
V\tilde{A}^{-1}\tilde{B}=C^{\sf T}
\end{align}
where  $\tilde{A} :=\tau^{-1}A$, and $\tilde{B} := \tau^{-1}B$.
If $G^{\rm red}$ is zero, then
\begin{align*}
k_{ij}(\delta_{ij}^{\star}) =
k_{ji}(\delta_{ji}^{\star})
,\quad
h_{ij}(\delta_{ij}^{\star}) = 
- h_{ji}(\delta_{ji}^{\star}),
\quad
h_{ii}(\delta_{ii}^{\star}) = 0,
\end{align*}
which means that $L$ is symmetric.
Furthermore, for
\[
\tilde{A} = 
%\mat{
%\sfdiag \left(  \frac{ X_i - X_i' }{ \tau_{{\rm d}i} } \right) & 0 \\
%0 & \sfdiag \left(  \frac{ X_i - X_i' }{ \tau_{{\rm q}i} } \right)
%}
\tau^{-1} \left\{ I_2 \otimes \sfdiag \left(  X_i -X_i'  \right) \right\}
\hat{A}
\]
where $X_i > X_i'$, we see that $\hat{A}$ is negative definite, because of the stability of $\hat{H}$.
Note that $\hat{A}$ is symmetric because
\[
\spliteq{
\tfrac{\partial g_{{\rm q}i} }{\partial E_{{\rm q}j} }(z^{\star})
&=
\tfrac{\partial g_{{\rm q}j} }{\partial E_{{\rm q}i} }(z^{\star})
= 
\tfrac{\partial g_{{\rm d}i} }{\partial E_{{\rm d}j} }(z^{\star})
=
\tfrac{\partial g_{{\rm d}j} }{\partial E_{{\rm d}i} }(z^{\star})
= 
k_{ij}(z^{\star}) \\
\tfrac{\partial g_{{\rm q}i} }{\partial E_{{\rm d}j} }(z^{\star})
&=
\tfrac{\partial g_{{\rm d}j} }{\partial E_{{\rm q}i} }(z^{\star})
= 
-h_{ij}(z^{\star})
}
\]
if $G^{\rm red}$ is zero.
Furthermore, $\hat{B}$ is equal to $-C^{\sf T}$ because 
\[
\spliteq{
\tfrac{\partial g_{{\rm q}i} }{\partial \delta_j }(z^{\star})
&=
- \tfrac{\partial P_j }{\partial E_{{\rm q}i} }(z^{\star})
= 
h_{ij}(\delta_{ij}^{\star})E_{{\rm q}j}^{\star}
+ k_{ij}(\delta_{ij}^{\star})E_{{\rm d}j}^{\star}
\\
\tfrac{\partial g_{{\rm q}i} }{\partial \delta_i }(z^{\star})
&=
- \tfrac{\partial P_i }{\partial E_{{\rm q}i} }(z^{\star})
= 
-\textstyle \sum\limits_{j\neq i}^N 
\bigl\{
h_{ij}(\delta_{ij}^{\star})E_{{\rm q}j}^{\star}
+ k_{ij}(\delta_{ij}^{\star})E_{{\rm d}j}^{\star}
\bigr\}
\\
\tfrac{\partial g_{{\rm d}i} }{\partial \delta_j }(z^{\star})
&=
- \tfrac{\partial P_j }{\partial E_{{\rm d}i} }(z^{\star})
= 
h_{ij}(\delta_{ij}^{\star})E_{{\rm d}j}^{\star}
- k_{ij}(\delta_{ij}^{\star})E_{{\rm q}j}^{\star}
\\
\tfrac{\partial g_{{\rm d}i} }{\partial \delta_i }(z^{\star})
&=
- \tfrac{\partial P_i }{\partial E_{{\rm d}i} }(z^{\star})
= 
-\textstyle \sum\limits_{j\neq i}^N 
\bigl\{
h_{ij}(\delta_{ij}^{\star})E_{{\rm d}j}^{\star}
- k_{ij}(\delta_{ij}^{\star})E_{{\rm q}j}^{\star}
\bigr\}.
}
\]
Therefore, $V$ can be chosen as the positive definite matrix $-\hat{A}$ satisfying \eqref{eq:cndQ}.
This proves the negative imaginaryness of $\hat{H}$.

Next, we consider $H$.
Because $\hat{H}$ is supposed to be stable, the pole of $H$ on the imaginary axis is only the origin, and it is simple.
Therefore, $H$ is positive real if and only if 
\begin{align}\label{eq:Gpr}
H(\bm{j}\omega) + H^{\sf T}(-\bm{j}\omega) \succeq 0
\end{align}
for all $\omega \in \mathbb{R}\setminus\{0\}$, and
\begin{align}\label{eq:cndG0}
\lim_{s\rightarrow 0} s H(s) = \lim_{s\rightarrow 0} \{ s H(s)\}^{\sf T}\succeq 0.
\end{align}
If $G^{\rm red}$ is zero, i.e., if $\hat{H}$ is negative imaginary, then
\begin{align}\label{eq:NIineq}
H(\bm{j}\omega) &+ H^{\sf T}(-\bm{j}\omega)  =
\frac{\bm{j}}{\omega} 
\left\{
\hat{H}(\bm{j}\omega) - \hat{H}^{\sf T}(-\bm{j}\omega)
\right\}
\end{align}
for all $\omega \in \mathbb{R}\setminus\{0\}$, which proves \eqref{eq:Gpr}.
Furthermore, we see that \eqref{eq:cndG0} is equivalent to \eqref{eq:pdsp} because
\begin{align*}
\lim_{s\rightarrow 0} s H(s) =
L - C\tilde{A}^{-1}\tilde{B} = L_0.
\end{align*}
Note that if the transmission network is lossless, both $L$ and
\begin{align*}
C\tilde{A}^{-1}\tilde{B} = C V^{-1} C^{\sf T}
\end{align*}
are symmetric.
Thus, the symmetry in \eqref{eq:cndG0} is also proven.

Finally, we prove that, if the transmission network is not lossless, or if \eqref{eq:pdsp} does not hold, then $H$ is not positive real.
The latter implication is trivial because \eqref{eq:pdsp} is equivalent to \eqref{eq:cndG0}.
For the former implication, if $G^{\rm red}$ is not zero, or equivalently, if $\hat{H}$ is not negative imaginary, then there exist some $\omega_0\geq 0$ and a sufficiently small $\epsilon >0$ such that
\begin{align*}
\lambda_{\rm min}\left[
\bm{j}
\left\{
H(\bm{j}(\omega_0 + \alpha )) - H^{\sf T}(-\bm{j}(\omega_0 + \alpha ))
\right\}
\right] < 0
\end{align*}
for any value $\alpha \in (0,\epsilon] $, 
where $\lambda_{\rm min}$ denotes the minimum eigenvalue.
Thus, \eqref{eq:cndQ} does not hold for $\omega \in (\omega_0, \omega_0+\epsilon] $.
\end{IEEEproof}

\vskip 1em

An important implication of Theorem~\ref{thm:EdynNI} is that a lossless transmission network is ``necessary" for the linearized version ${\mathds G}_{\vartriangle}$ in \eqref{eq:Gss} to be passive.
This also means that the nonlinear original $\mathds{G}$ in \eqref{eq:sys2all} is never EI passive if the transmission is lossy.
In fact, the positive semidefiniteness of $L_0$ in  \eqref{eq:pdsp}, which is one of the necessary conditions for ${\mathds G}_{\vartriangle}$ to be passive, can be understood as the ``convexity" of the strain energy function $U$ in  \eqref{eq:defstEkron}.
We will discuss this point in the next section.

\section{A New Link Between Two-Axis and Classical Generator Models in Terms of EI Passivity}\label{sec:conv}

In this section, we will prove Facts (b) and (c) in Section~\ref{sec:intro}, namely the convex domain of the strain energy function is equal to the ``largest" set of equilibria over which the electromagnetic subsystem is EI passive, and the EI passivity of the two-axis model network is equivalent to that of the classical model network derived by the SPA, and the stability of its flux linkage dynamics.
Throughout this section, we again assume that the transmission network is lossless because it is necessary for the EI passivity.

\subsection{Largest Set of Equilibria for EI Passivity}

First, we consider Fact (b) based on the analysis of the linearized electromagnetic subsystem in Section~\ref{sec:anakronlin}.
Recall that the non-negativity of the storage function $W_{z^*}$ in \eqref{eq:stfunW} is equivalent to the convexity of the strain energy function $U$.
It is interesting to note that the Hessian of $U$ is found to be
\begin{equation}\label{eq:HessU}
\nabla^2 U(z^{\star})
=
\mat{
L & -\hat{B}^{\sf T} \\
-\hat{B} & -\hat{A}
}
\end{equation}
where $\hat{A}$, $\hat{B}$, and $L$ are system matrices of the linearized version defined as in \eqref{eq:sysmats}.
From this coincidence, the following important fact can be deduced.

\begin{theorem}\label{thm:Hesconv2}
The equilibrium set $\mathcal{E}$ in \eqref{eq:eqEG}
is the largest set of equilibria over which the electromagnetic subsystem $\mathds{G}$ in \eqref{eq:sys2all} is EI passive.
\end{theorem}

\begin{IEEEproof}
To prove the claim, it suffices to show that, for any feasible equilibrium $z^{\star}$ not belonging to $\mathcal{E}$, $H$ in \eqref{eq:trGs}, which is an implicit function of $z^{\star}$, is not positive real.
This is equivalent to prove that $\nabla^2 U(z^{\star})$ is positive semidefinite for any feasible equilibrium $z^{\star}$ such that $H$ is positive real.
This is also equivalent to prove that, if $H$ is positive real, then $L_0$ is positive semidefinite, and $\hat{A}$ is negative definite.

As shown in Theorem~\ref{thm:EdynNI}, the positive semidefiniteness of $L_0$ is necessary for $H$ to be positive real.
In addition, the negative definiteness of $\hat{A}$ is also necessary because, if $\hat{A}$ has a positive eigenvalue, then $H$ is not stable, and, if $\hat{A}$ has a zero eigenvalue, then the pole of $H$ at the origin is not simple.
Therefore, considering the Schur complement of \eqref{eq:HessU} with respect to $-\hat{A}$, we can see that, if $H$ is positive real, then $\nabla^2 U(z^{\star})$ is positive semidefinite.
\end{IEEEproof}

Theorem~\ref{thm:Hesconv2} gives a remarkable link between the nonlinear analysis and linear analysis.
In fact, the convex domain of the strain energy function is found to be the ``largest" set of equilibria whose stability can be proven by virtue of the EI passivity.
It should be emphasized that such a largest equilibrium set is generally difficult to find by nonlinear analysis because searching all possible storage functions is not realistic.

\subsection{New Link Between Two-Axis and Classical Models}

%\begin{figure*}[ht]
%\centering
%\includegraphics[width = .85\linewidth]{ sysdiagram.pdf}
%%\medskip
%\caption{Relations among electromagnetic subsystems, strain energy functions and storage functions. 
%``Kron" indicates the Kron reduction, ``Lin" indicates the linearization, and ``SPA" indicates the singular perturbation approximation of the flux linkage dynamics such that $(\tau_{\rm q},\tau_{\rm d}) \rightarrow 0$.}
%\label{figsysdiagram}
%\end{figure*}

Next, we will prove Fact (c).
In particular, we will show that the positive semidefiniteness of $L_0$ in \eqref{eq:Lapmat} can be understood as the convexity of the strain energy function of the ``classical" model network.

Suppose that the transmission network is lossless.
Then, the electromagnetic subsystem of the classical model network is obtained as the collection of
\begin{equation}\label{eq:Gtilred}
\tilde{{\mathds G}}_i :
\simode{
\dot{\delta}_i&= \tilde{u}_i \\
\tilde{y}_i & \textstyle = - \sum\limits_{j=1}^N V_{{\rm fd}i}^{\star}V_{{\rm fd}j}^{\star} \tilde{B}^{\rm red}_{ij} \sfsin \delta_{ij} 
}
\end{equation}
where $\tilde{B}^{\rm red}_{ij}$ is the $(i,j)$-element of the reduced susceptance matrix defined with the synchronous reactances as
\[
\tilde{B}^{\rm red} := \sfIm 
\Bigl[
- \bm{j}
\bigl\{
%\underbrace{
\sfdiag \left( X_i \right) 
-  \bm{j} \sfdiag \left( X_i \right) \overline{\bm{Y}} \sfdiag \left( X_i \right)
%}_{\tilde{\bm{\varGamma}}}
\bigr\}^{-1}
\Bigr].
\]
Note that  $\tilde{B}^{\rm red}$ is not identical to $B^{\rm red}$ in \eqref{eq:GBred} as the synchronous reactance $X_i$ is different from the transient reactance $X_i'$.
Then, the strain energy function of $\tilde{{\mathds G}}$ is found as
\begin{equation}\label{eq:Udel}
\tilde{U} (\delta) :=  \frac{1}{2} \sum_{i=1}^N \sum_{j=1}^N V_{{\rm fd}i}^{\star}V_{{\rm fd}j}^{\star} \tilde{B}_{ij}^{\rm red} \sfcos \delta_{ij}.
\end{equation}
We can verify that $\tilde{{\mathds G}}$ is EI passive for the storage function
\[
\tilde{W}_{\delta^{\star}} ( \delta ) :=
\mathfrak{B}_{\delta^{\star}} \bigl[
\tilde{U} ( \delta ) 
\bigr]
\]
with respect to the equilibrium set
\begin{equation}\label{eq:calEcl}
\tilde{\mathcal{E}} :=
{\sf int} \! \left\{
\delta^{\star} \in \mathbb{S}^N:
\nabla^2 \tilde{U} (\delta^{\star}) \succeq 0
\right\}.
\end{equation}
On the top of the existing relationship between the two-axis model and the classical model in Proposition~\ref{prop:2axcl}, we can prove the following novel fact.

\begin{theorem}\label{thm:Hesconcl}
Consider the electromagnetic subsystem $\mathds{G}$ in \eqref{eq:sys2all}, and the system matrix $\hat{A}$ in \eqref{eq:sysmats}.
Then, the equilibrium set $\mathcal{E}$ in \eqref{eq:eqEG} is identical to
\begin{equation}\label{eq:eqEGredT}
\mathcal{E} =
{\sf int} \! \left\{
z^{\star}
\in \mathcal{Z}^{\star} :
\hat{A}(z^{\star}) \preceq 0,\ 
\nabla^2 \tilde{U}( \delta^{\star}) \succeq 0
\right\}
\end{equation}
where $\tilde{U}$ in \eqref{eq:Udel} corresponds to the strain energy function of the electromagnetic subsystem $\tilde{{\mathds G}}$ in \eqref{eq:Gtilred}.
\end{theorem}

\begin{IEEEproof}
Considering the Schur complement of \eqref{eq:HessU} with respect to $-\hat{A}$, we have
\[
\mathcal{E} =
{\sf int} \! \left\{
z^{\star}
\in \mathcal{Z}^{\star} :
\hat{A}(z^{\star}) \preceq 0,\ 
L_0(z^{\star}) \succeq 0
\right\}
\]
where $L_0$ is defined as in \eqref{eq:Lapmat}, which is equal to
\[
L_0(z^{\star}) = L(z^{\star})+\hat{B}^{\sf T}(z^{\star})
\hat{A}^{-1}(z^{\star})\hat{B}(z^{\star}).
\]
Therefore, it suffices to show that, if $z^{\star}\in \mathcal{Z}^{\star} $, then
\begin{equation}\label{eq:L0nabU}
L_0(z^{\star}) = \nabla^2 \tilde{U}( \delta^{\star}).
\end{equation}
To prove this, we consider the limit such that $\tau$ of ${\mathds G}_{\vartriangle}$ in \eqref{eq:Gss} is sufficiently small.
Then, we have 
\begin{subequations}\label{eq:Gsstillin}
\begin{equation}
\tilde{{\mathds G}}_{\vartriangle}: \simode{
\dot{{\delta}}_{\vartriangle}
&= \textstyle 
\tilde{u}_{\vartriangle}  \\
\tilde{y}_{\vartriangle} &= L_0( z^{\star}) \delta_{\vartriangle}. 
}
\end{equation}
From Proposition~\ref{prop:2axcl}, we see that this $\tilde{{\mathds G}}_{\vartriangle}$ corresponds to the linearized version of the electromagnetic subsystem of the classical model network, which can also be derived by the linearization of $\tilde{{\mathds G}}$ in \eqref{eq:Gtilred} as
\begin{equation}
\tilde{{\mathds G}}_{\vartriangle}: \simode{
\dot{{\delta}}_{\vartriangle}
&= \textstyle 
\tilde{u}_{\vartriangle}  \\
\tilde{y}_{\vartriangle} &= \nabla^2 \tilde{U}( \delta^{\star}) \delta_{\vartriangle}. 
}
\end{equation}
\end{subequations}
Thus, the equality in \eqref{eq:L0nabU} is proven.
\end{IEEEproof}

Theorem~\ref{thm:Hesconcl} states that the EI passivity of the electromagnetic subsystem of the two-axis model network can be characterized by two conditions.
One is the stability of $\hat{A}$, which is the Jacobian of the vector field of the flux linkage dynamics, and the other is the convexity of $\tilde{U}$, which is the strain energy function of the ``classical" model network derived by the SPA of the flux linkage dynamics.
It should be emphasized again that the reduced susceptance $B^{\rm red}_{ij}$ in \eqref{eq:stfunW} is different from the reduced susceptance $\tilde{B}^{\rm red}_{ij}$ in \eqref{eq:Udel}.
Nevertheless, it is proven that the equilibrium set $\mathcal{E}$ in \eqref{eq:eqEG} is identical to that in \eqref{eq:eqEGredT}.

%The relations among the different expressions are summarized in Fig.~\ref{figsysdiagram}.

%To prove \eqref{eq:L0nabU}, we have used the fact that $\tilde{{\mathds G}}^{\rm red}_{\vartriangle}$ in \eqref{eq:Gsstillin} can be derived by two ways.
%For reference, the strain energy function for the DAE representation of the salient pole type classical generator model is given as
%\[
%\spliteq{
%\tilde{U} (\delta, |\bm{V}|, \angle \bm{V} ) &:=
% \sum_{i=1}^N 
%\biggl\{
% - 
%\frac{1}{4}
%\left(
%\frac{1}{X_{{\rm q}i}}
%-
%\frac{1}{X_{{\rm d}i}}
%\right)
%\bigl( \bm{V}_{{\rm q}i}^2 - \bm{V}_{{\rm d}i}^2 \bigr)\\
%& 
%- \frac{V_{{\rm fd}i}^{\star} }{X_{{\rm d}i}} \bm{V}_{{\rm q}i}
%+ 
%\frac{1}{4}
%\left(
%\frac{1}{X_{{\rm q}i}}
%+
%\frac{1}{X_{{\rm d}i}}
%\right) |\bm{V}_i|^2
%\\
%&- 
%\sum_{j=1}^{N} 
%\frac{B_{ij} }{2} |\bm{V}_i| |\bm{V}_j| \sfcos(\angle \bm{V}_i -\angle \bm{V}_j)
%\biggr\}.
%}
%\]

It is interesting to note that $L_0( \delta^{\star})$ in \eqref{eq:Lapmat}, or equivalently $\nabla^2 \tilde{U}( \delta^{\star})$ in \eqref{eq:eqEGredT}, is a weighted graph Laplacian corresponding to the ``spring stiffness matrix" of the linearized version of the classical model network.
In particular, it is written as
\[
M \ddot{\delta}_{\vartriangle} 
+ D \dot{\delta}_{\vartriangle} + 
\omega_0 L_0( \delta^{\star}) \delta_{\vartriangle} = 0,
\]
where $M$ and $D$ are the diagonal matrices composed of $M_i$ and $D_i$, respectively.
This means that the convexity of the strain energy function of the classical model network is equivalent to the positive semidefiniteness of the spring stiffness matrix in the linear second-order system.

Furthermore, $L_0( \delta^{\star})$ can be viewed as a matrix generalization of the \textit{synchronizing torque coefficients} \cite{sauer2017power}.
This is because its $(i,j)$-element is equal to
$\tfrac{\partial \tilde{P}_i}{\partial \delta_j}( \delta^{\star})$ where 
\[
\tilde{P}_i ( \delta)= 
- \sum_{j=1}^N V_{{\rm fd}i}^{\star}V_{{\rm fd}j}^{\star} \tilde{B}^{\rm red}_{ij} \sfsin \delta_{ij}
\]
denotes the active power output of $\tilde{{\mathds G}}$ in \eqref{eq:Gtilred}.
Note that the conventional concept of synchronizing torque coefficients is well defined for the single-machine infinite-bus model with the classical generator model.
It is interesting to note that such a basic stability concept for the single classical model can be generalized by  analyzing the EI passivity for the two-axis model network.
Readers interested in a simple sufficient condition for the matrix version $L_0( \delta^{\star})$ or $\nabla^2 \tilde{U}( \delta^{\star})$ to be positive semidefinite are referred to Appendix~\ref{sec:EIclass}.

For the stability analysis of lossy power systems, the convexity of the strain energy function should be replaced by the positive semidefiniteness of the synchronizing torque coefficient matrix $L_0( \delta^{\star})$, which can be calculated even in the lossy case.
We will demonstrate such a generalization in Section~\ref{sec:numex} through a numerical simulation.

%\begin{remark}
%The storage functions of the linearized ${\mathds G}_{\vartriangle}$ and $\tilde{\mathds G}_{\vartriangle}$ can be found as the quadratic functions of
%\[
%\spliteq{
%W_{\vartriangle}(z_{\vartriangle})
%& :=
%z_{\vartriangle}^{\sf T}
%\nabla^2 U(z^{\star})
%z_{\vartriangle},
%\\
%\tilde{W}_{\vartriangle}(\delta_{\vartriangle})
%&:=
%\delta_{\vartriangle}^{\sf T}
%\nabla^2 \tilde{U} (\delta^{\star})
%\delta_{\vartriangle}
%}
%\]
%where $z_{\vartriangle}$ denotes the vector composed of $\delta_{\vartriangle}$ and $E_{\vartriangle}$.
%\end{remark}

\section{Numerical Simulation}\label{sec:numex}

\begin{figure}[t]
\centering
\includegraphics[width = .75\linewidth]{ 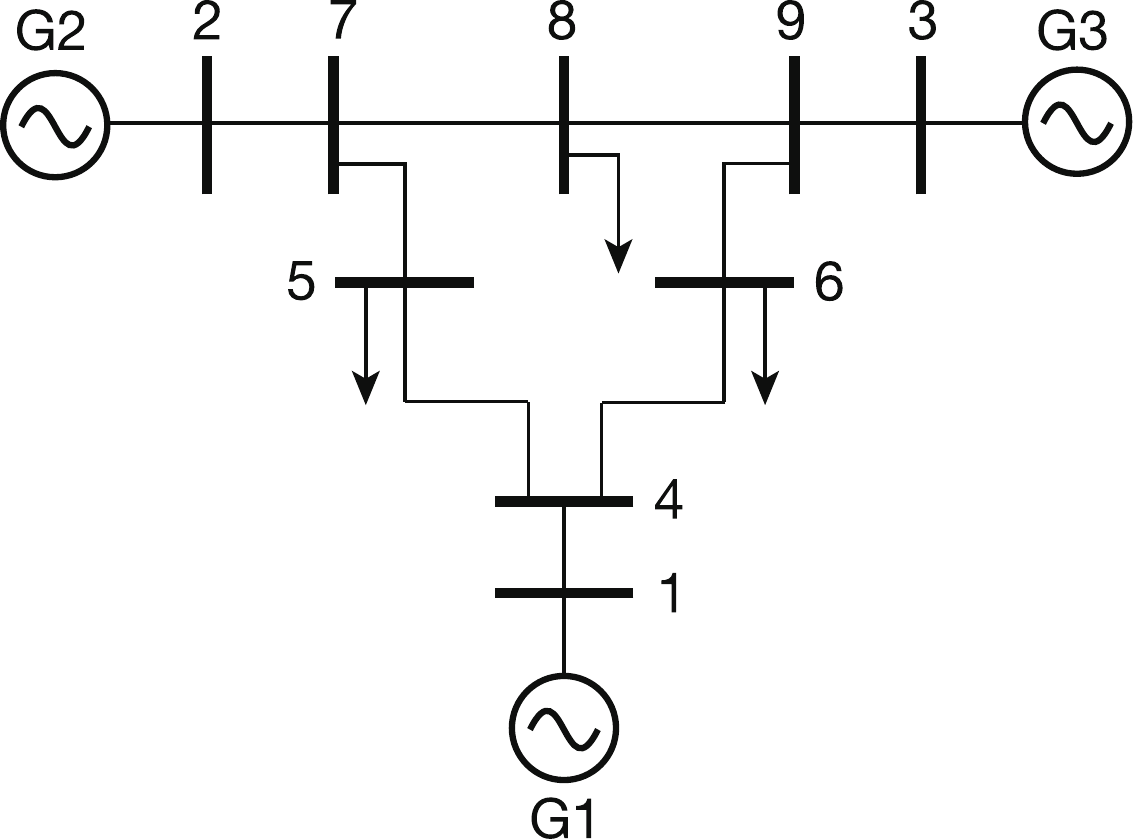}
%\medskip
\caption{IEEE 9-bus system model. 
}
\label{figWECC9bus}
\end{figure}

In this section, we demonstrate practical significance of our mathematical analysis using the IEEE 9-bus system model shown in Fig.~\ref{figWECC9bus}.
The values of the constant model parameters associated with the generators, loads, and transmission lines are listed in Tabs.~\ref{table:genconst} and \ref{table:transconst}.
The generator and transmission line constants are taken from Example~7.1 of \cite{sauer2017power}.
Note that Buses~4, 7, and 9 can be removed in an equivalent manner \cite{ishizaki2018graph} so that a generator or load is connected to each of all buses in the resulting power system.

The generators are supposed to be the two-axis model.
Two cases are considered for the load models.
One is the case of induction motor loads represented by the classical model in \eqref{eq:clmodel}, whose constants are taken from \cite{chen2020modelling} for virtual synchronous generators.
The other is the case of inverter loads with power-frequency droop control, represented by
\begin{equation}
D_i \dot{\delta}_i = \omega_0 (P_{{\rm m}i}-P_i)
\end{equation}
which is obtained as a special case of the classical model where the inertia constant $M_i$ is sufficiently small.
For the loads, the active power references are constants specified as
\[
P_{{\rm m}5}=-1.25,\quad  
P_{{\rm m}6}=-0.90,\quad
P_{{\rm m}8}=-1.00.
\]

For the transmission network, we consider both lossless and lossy cases.
In the lossless case, we assume that the conductance of each transmission line is zero.
The objective of this numerical simulation is to demonstrate that the largest set of equilibria over which the electromagnetic subsystem is EI passive, or equivalently the convex domain of its strain energy function, is almost identical to the set of all stable equilibria, i.e., Fact (d) in Section~\ref{sec:intro}.
In addition, this can be generalized to the lossy power system.

\begin{table}[b]
%\medskip
 \caption{Generator and load constants.}
 \label{table:genconst}
 \centering
  \begin{tabular}{lrrrrrrrrrrrrr}
   \hline
 & Bus~1 & Bus~2 & Bus~3  & Bus~5 & Bus~6 & Bus~8 \\   \hline \hline%& Bus~7 & Bus~8 & Bus~9  \\   \hline \hline
$M_i$ [s] &  0.1254 & 0.0340  & 0.0160 & 0.0042 & 0.0042 & 0.0042 \\
$D_i$ [pu] & 0.0125 & 0.0068 & 0.0048  & 0.0003 & 0.0003 & 0.0003 \\
$X_{i}$ [pu] & 0.1460 & 0.8958 & 1.3120 & 0.3 & 0.3 & 0.3\\
$X_{i}'$ [pu] & 0.0608 & 0.1198 & 0.1813 \\
%$X_{{\rm q}i}$ [pu] &0.0969 & 0.8645 & 1.2578 \\ 
%$X_{{\rm q}i}'$ [pu] & 0.0969 & 0.1969 & 0.2500 \\
$\tau_{{\rm d}i}$ [s] & 8.9600 & 6.0000 & 5.8900 \\
$\tau_{{\rm q}i}$ [s] & 0.3100 & 0.5350 & 0.6000 \\
%$P_{{\rm ld}i}$ [pu] &&&&& $-$1.2500 & $-$0.9000  && $-$1.0000 \\
%$Q_{{\rm ld}i}$ [pu] &&&&& $-$0.5000  & $-$0.3000 && $-$0.3500 \\   
%$\beta_i$ [pu] & 0.0000 & 0.0000 & 0.0000 & 0.1670 & 0.2410 & 0.2580 & 0.2270 & 0.1790 & 0.2830 \\
\hline
  \end{tabular}
\end{table}

\begin{table*}[t]
%\medskip
 \caption{Transmission line constants.}
 \label{table:transconst}
 \centering
  \begin{tabular}{lrrrrrrrrrrrrrrrrl}
   \hline
 & $(1,4)$ & $(2,7)$ & $(3,9)$ & $(4,5)$ & $(4,6)$ & $(5,7)$ & $(6,9)$ & $(7,8)$ & $(8,9)$ & \\   \hline \hline
$g_{ij}$ [pu] & 0.0000  & 0.0000 & 0.0000 & 1.3650 & 1.9420 & 1.1880 & 1.2820 & 1.6170 & 1.1550 & \\
$b_{ij}$ [pu] & $-$17.361 & $-$16.000 & $-$17.065 & $-$11.604 & $-$10.511 & $-$5.9750 & $-$5.5880 & $-$13.698 & $-$9.7840 & \\
$c_{ij}$ [pu] &  0&  0&  0& 0.4669&  0.4191&  0.8117&  0.9496&  0.3952&  0.5544 & {\footnotesize$\times10^{-3}$}\\  \hline
  \end{tabular}
\end{table*}

\begin{figure*}[t]
\centering
\includegraphics[width = .99\linewidth]{ 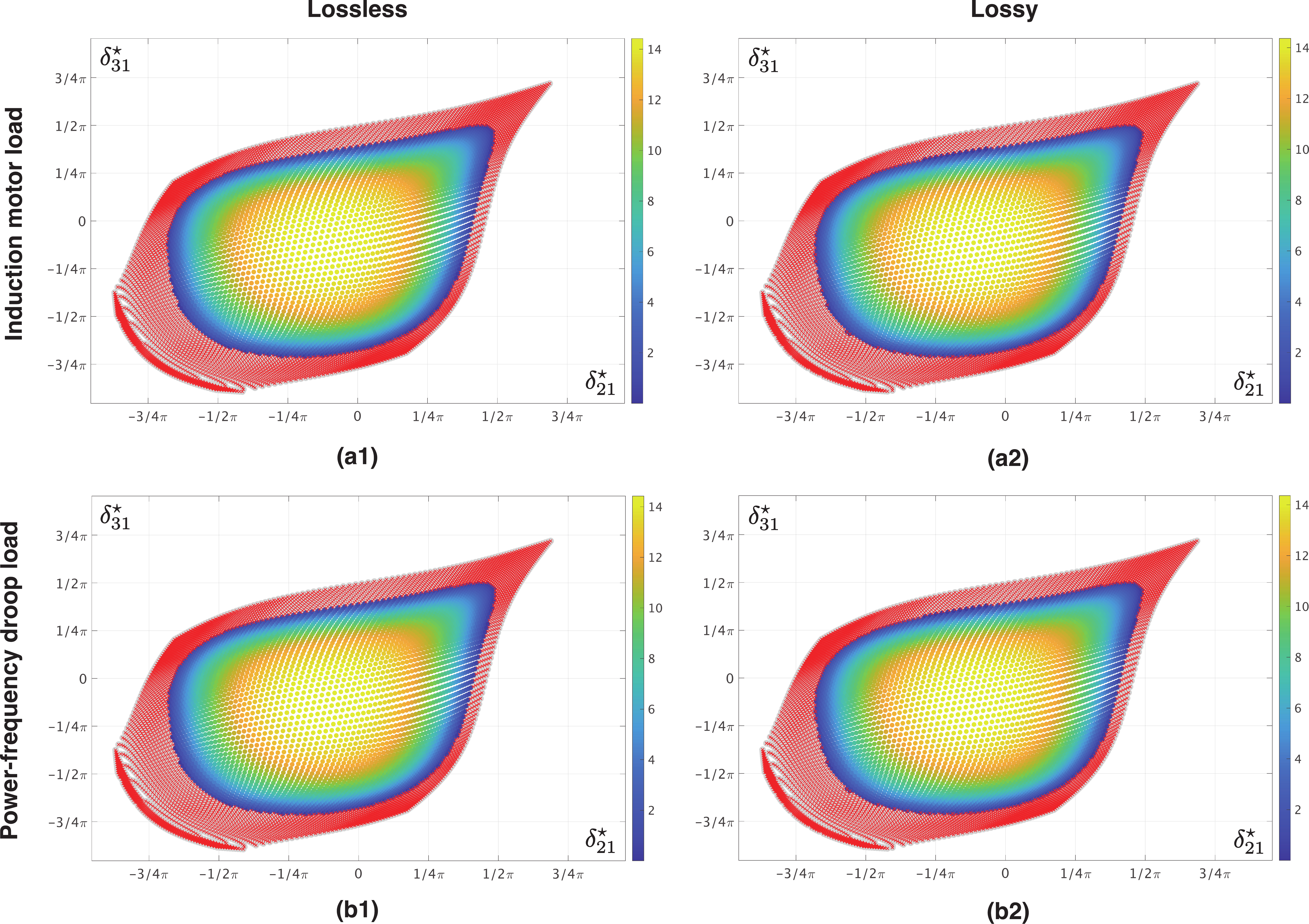}
\caption{Analysis of stable equilibria for lossless and lossy power systems with the two-axis generator model. 
(a1) Lossless power system with induction motor load model.
(b1) Lossless power system with power-frequency droop load model.
(a2) Lossy power system with induction motor load model.
(b2) Lossy power system with power-frequency droop load model.
}
\label{DomainAll}
\end{figure*}

\subsection{The Lossless Case}

We numerically study the stability of all feasible equilibria.
First, we consider the case of the lossless power system.
The results are shown in Fig.~\ref{DomainAll}(a1) for the induction motor load model and Fig.~\ref{DomainAll}(b1) for the power-frequency droop load model, where the plane of $(\delta_{21}^{\star},\delta_{31}^{\star})$ is gridded.
The colored area represents the feasible equilibrium set of the entire power system model, the area of the red crosses represents the unstable subset, and the area of the blue to yellow dots represents the equilibrium set $\mathcal{E}$ in \eqref{eq:eqEG} or equivalently \eqref{eq:eqEGredT}.
The color bar represents the products of all positive eigenvalues of $L_0(\delta^{\star})$ in in \eqref{eq:Lapmat}, which quantifies the degree of synchronizing torque.
We can see that the equilibrium set $\mathcal{E}$ exactly captures the set of all stable equilibria without any  exception for both load models.
In fact, the results in Figs.~\ref{DomainAll}(a1) and (b1) are exactly the same.

\subsection{The Lossy Case}

Next, we analyze the stability of the lossy power system.
Using the Jacobian of the flux linkage dynamics and the synchronizing torque coefficient matrix, i.e., $A(z^{\star})$ in \eqref{eq:sysmats} and  $L_0( \delta^{\star})$ in \eqref{eq:Lapmat}, we define the equilibrium set
\begin{equation}\label{eq:Eplus}
\mathcal{E}_+ :=
{\sf int} \! \left\{
z^{\star}
 \in \mathcal{Z}^{\star} :
\mbox{\!\!
$A(z^{\star})$ and $-L_0(\delta^{\star})$ are stable
\!}
\right\}
\end{equation}
where the stability means that the eigenvalues lie in the closed left half-plane.
This $\mathcal{E}_+$ is a generalization of $\mathcal{E}$, which can be computed even in the lossy case.

The results are shown in Fig.~\ref{DomainAll}(a2) for the induction motor load model and Fig.~\ref{DomainAll}(b2) for the power-frequency droop load model.
As can be seen from these results, $\mathcal{E}_+$ in \eqref{eq:Eplus} almost exactly captures the set of all stable equilibria, with few exceptions around the dark blue boundary.
These results demonstrate that the convex domain of the strain energy function, or equivalently the domain such that the synchronizing torque coefficient matrix is positive semidefinite, is almost identical to the set of all stable equilibria even in the case of lossy power systems.
This is not trivial, because the set of equilibria over which a system is EI passive is in general only a subset of the stable equilibria.

\section{Concluding Remarks}\label{sec:conc}

In this paper, we have analyzed the EI passivity of a multi-machine power system considering two-axis generator models.
The key practical findings are summarized as follows.
\begin{itemize}
\item The convexity of a strain energy function characterizes the largest set of equilibria over which an electromagnetic subsystem is EI passive, which is almost identical to the set of all stable equilibria.
\item The convexity of the strain energy function is equivalent to the positive semidefiniteness of a synchronizing torque coefficient matrix defined for the two-axis model network, and the stability of the flux linkage dynamics.
\item The positive semidefiniteness of the synchronizing torque coefficient matrix characterizes the set of all stable equilibria even for lossy power systems, although the EI passivity does not hold in the case of lossy transmission.
\end{itemize}
These practical findings support the rationality of discussing the stability of even lossy power systems in terms of the EI passivity.
Future challenges include the application of this stability criterion to the calculation of optimal power flow and the design of smart inverters for systems with heterogeneous components such as renewable energy.

%\appendix

\appendices

\section{Mathematical Details of Kron Reduction}\label{sec:kron}

\subsection{Derivation of Equivalent ODE Representation}

We explain the details of the derivation of the ODE equivalent in Section~\ref{sec:feedrep} from the original DAE representation in Section~\ref{sec:DAEmod}.
This procedure is often called the \textit{Kron reduction} for phasor circuits, where static nodes or buses in a power system are removed in a mathematically equivalent procedure when the generators are considered as dynamical nodes.
This corresponds to a phasor circuit version of the standard $Y$-$\Delta$ transformation in circuit theory \cite{ishizaki2018graph,kron1939tensor}.

The goal here is to eliminate the bus current and voltage phasors $\bm{I}_i$ and $\bm{V}_i$ involved in \eqref{eq:genmodel}.
From \eqref{eq:phVI}, we see that
\[
\bm{I}_{\rm q} + \bm{j} \bm{I}_{\rm d} =
-\bm{j} \sfdiag \left(
\tfrac{1}{X_i'}
\right)
\left\{
\bm{V}_{\rm q} + \bm{j} \bm{V}_{\rm q}
-
\bigl(
E_{\rm q} + \bm{j} E_{\rm d} 
\bigr)
\right\}.
\]
Furthermore, from the network equation in \eqref{eq:ohmY}, we have
\[
\bm{I}_{\rm q} + \bm{j} \bm{I}_{\rm d} =
\sfdiag \left(
e^{\bm{j}\delta_i}
\right)
\overline{\bm{Y}}
\sfdiag \left(
e^{-\bm{j}\delta_i}
\right)
\bigl(
\bm{V}_{\rm q} + \bm{j} \bm{V}_{\rm d}
\bigr).
\]
Eliminating $\bm{I}_{\rm q}$ and $\bm{I}_{\rm d}$ from them, we have
\begin{equation}\label{eq:kronGam}
\spliteq{
\sfdiag \left(
\tfrac{1}{X_i'}
\right)
&
\bigl(
\bm{V}_{\rm q} + \bm{j} \bm{V}_{\rm d}
\bigr)
\\
& =
\underbrace{
\sfdiag\left(
e^{\bm{j} \delta_i}
\right)
\bm{\varGamma}^{-1}
\sfdiag\left(
e^{-\bm{j} \delta_i}
\right)
}_{\star}
\bigl(
E_{\rm q} + \bm{j} E_{\rm d} 
\bigr)
}
\end{equation}
where $\bm{\varGamma}$ is defined as 
\begin{equation}\label{eq:dGamma}
\bm{\varGamma}:=
\sfdiag \left( X_i' \right) 
-  \bm{j} \sfdiag \left( X_i' \right) \overline{\bm{Y}} \sfdiag \left( X_i' \right).
\end{equation}
This $\bm{\varGamma}$ is non-singular in a standard parameter setting, as will be proven in Appendix~\ref{sec:supmath}.
Note that the reduced admittance matrix $\bm{Y}^{\rm red}$ in \eqref{eq:Yred} is equal to $-\bm{j} \bm{\varGamma}^{-1} $.

Let $\bm{\gamma}_{ij}^{-1}$ denote the $(i,j)$-element of $\bm{\varGamma}^{-1}$.
Then, the $(i,j)$-element of ``$\star$" in \eqref{eq:kronGam}
is written as $\bm{\gamma}_{ij}^{-1} e^{\bm{j}\delta_{ij}} $.
Thus, its real and imaginary parts are, respectively, found to be $k_{ij}$ and $h_{ij}$ in \eqref{eq:defkh} where 
$B_{ij}^{\rm red}$ and $G_{ij}^{\rm red}$ are equal to
\[
B_{ij}^{\rm red} = 
- \sfRe \left[ \bm{\gamma}_{ij}^{-1}  \right]
,\quad
G_{ij}^{\rm red} = 
\sfIm \left[
\bm{\gamma}_{ij}^{-1} 
\right].
\]
Thus, we can rewrite the real part of \eqref{eq:kronGam} as
\[
\frac{\bm{V}_{{\rm q}i} }{X_i'} =
\sfRe \left[ \textstyle
\sum_{j=1}^N\limits \bm{\gamma}_{ij}^{-1} e^{\bm{j}\delta_{ij}}  
\bigl(
E_{{\rm q}j} + \bm{j} E_{{\rm d}j} 
\bigr)
\right],
\]
and the imaginary part as
\[
\frac{\bm{V}_{{\rm d}i} }{X_i'} =
\sfIm \left[ \textstyle
\sum_{j=1}^N\limits \bm{\gamma}_{ij}^{-1} e^{\bm{j}\delta_{ij}}  
\bigl(
E_{{\rm q}j} + \bm{j} E_{{\rm d}j} 
\bigr)
\right] ,
\]
which are equal to $g_{{\rm q}i}$ and $g_{{\rm d}i}$ in \eqref{eq:kronterms}.
Finally, substituting these into \eqref{eq:genmodel}, we obtain the ODE equivalent in \eqref{eq:kron2ax}.

\subsection{Supplemental Materials}\label{sec:supmath}

We prove several mathematical properties on the conductance matrix $G$ and susceptance matrix $B$, i.e., the real and imaginary parts of the admittance matrix $\bm{Y}$.
The conductance matrix is a weighted graph Laplacian given as
\[
G_{ij}
=
\left\{
\begin{array}{cl}
-g_{ij}, & i\neq j \\
\sum\limits_{j=1}^{N} g_{ij},  & i=j
\end{array}
\right.
\]
where $G_{ij}$ denotes the $(i,j)$-element of $G$, and $g_{ij} $ denotes the line conductance being non-negative.
In a similar way, the susceptance matrix $B$ is  given as
\[
B_{ij}
=
\left\{
\begin{array}{cl}
b_{ij}, & i\neq j \\
\sum\limits_{j=1}^{N} \left( b_{ij} +  \frac{\omega_0c_{ij}}{2} \right), & i=j
\end{array}
\right.
\]
where $B_{ij}$ denotes the $(i,j)$-element of $B$, and $b_{ij} $ denotes the line susceptance being non-positive.
We can decompose $B$ into the sum of the sign inversion of a weighted graph Laplacian and a positive semidefinite diagonal matrix as
\begin{equation}\label{eq:defbeta}
B = B_0 + \sfdiag \biggl( 
\underbrace{
\textstyle
\sum\limits_{j=1}^N
\tfrac{\omega_0c_{ij}}{2}
}_{\beta_{i}}
\biggr)_{i \in \{1,\ldots,N\}
}.
\end{equation}
Without loss of generality, we assume that the transmission network is connected, i.e., a path of transmission lines exists between any pair of two buses.
This assures that the kernel of $B_0$ is one-dimensional. 
By definition, it is clear that 
\begin{equation}\label{eq:kerGB0}
\sfspan\{\mathds{1}\}
\subseteq 
\sfker G
,\quad
\sfspan\{\mathds{1}\}
=
\sfker B_0.
\end{equation}
Furthermore, $G$ is positive semidefinite, and $B_0$ is negative semidefinite.

As shown in the following lemma, $\bm{\varGamma}$ in \eqref{eq:dGamma} is non-singular in a proper parameter setting.

\begin{lemma}\label{lem:gamnon}
Consider the admittance matrix $\bm{Y}$ in \eqref{eq:GandB}. 
For the constants $\beta_1,\ldots,\beta_N$ in \eqref{eq:defbeta}, if
\begin{align}\label{eq:Xbeta}
\beta_i X'_i \leq 1
,\quad \forall i \in \{1,\ldots,N\}
\end{align}
and the strict inequality holds for at least one bus, then $\bm{\varGamma}$ in \eqref{eq:dGamma} is non-singular.
\end{lemma}

\begin{IEEEproof}
The real part and imaginary part of $\bm{\varGamma}$, denoted by $M$ and $N$, are written as
\begin{equation}\label{eq:defMNm}
\spliteq{
M &:= \sfdiag \bigl( X_i' (1- \beta_i X_i' ) \bigr) 
- \sfdiag \left( X_i' \right) B_0 \sfdiag \left( X_i' \right), \\
N &:= - \sfdiag \left( X_i' \right) G \sfdiag \left( X_i' \right).
}
\end{equation}
Because $B_0$ in \eqref{eq:defbeta} is negative semidefinite, if \eqref{eq:Xbeta} holds, then $M$ is positive semidefinite.
In particular, if 
\[
\beta_i X_i' <1
\]
for at least one $i\in \{1,\ldots,N\}$, then $M$ is positive definite because
\[
\underbrace{
\sfker \sfdiag \left( X_i' \right) B_0 \sfdiag \left( X_i' \right)
}_{\sfspan \left\{ \sfdiag \left( 1/X_i' \right) \mathds{1} \right\}}
\nsubseteq
\sfker \sfdiag \bigl( X_i' (1- \beta_i X_i' ) \bigr).
\]
Thus, from the symmetry of $N$, we see that 
\[
L_1 := M+NM^{-1}N
\]
is positive definite, meaning that $L_1$ is non-singular.
Furthermore, $M+\bm{j}N$ is non-singular if and only if 
\[
L_2 := \mat{M & -N \\ N & M}
\]
is non-singular.
Because the determinant of $L_2$ is 
\[
\sfdet (L_2) = \sfdet(M) \sfdet (L_1)
\]
and $\sfdet(M)$ is non-zero, $L_2$ is non-singular if and only if $L_1$ is non-singular.
This proves the claim.
\end{IEEEproof}

In a standard parameter setting, the ground capacitances are sufficiently small that \eqref{eq:Xbeta} generally holds.
As shown in the following proposition, these reduced matrices have the same definiteness of the original ones before the Kron reduction.

\begin{proposition}\label{prop:defBGred}
Consider the reduced admittance matrix $\bm{Y}^{\rm red}$ in \eqref{eq:Yred}.
If the inequality condition in Lemma~\ref{lem:gamnon} is satisfied, then the reduced conductance matrix  $G^{\rm red}$ is positive semidefinite, and the reduced susceptance matrix $B^{\rm red}$ is negative definite.
\end{proposition}

\begin{IEEEproof}
With the same notation as that in the proof of Lemma~\ref{lem:gamnon}, denote the real and imaginary parts of $\bm{\varGamma}$ by $M$ and $N$.
Furthermore, denote the real and imaginary parts of $\bm{\varGamma}^{-1}$ by $P$ and $Q$.
Note that 
\[
\bm{\varGamma}
\bm{\varGamma}^{-1} =I
\quad
\Longleftrightarrow
\quad
\underbrace{
\mat{
M & -N\\
N & M
}
}_{L_2}
\mat{
P & -Q\\
Q & P
}
=I,
\]
which means that the diagonal and off-diagonal blocks of $L_2^{-1}$ correspond to $P$ and $Q$, respectively.
In particular, we see that
\begin{align*}
L^{-1}_2=
\mat{
L_1^{-1} & M^{-1}NL_1^{-1} \\
-L_1^{-1}NM^{-1} & L_1^{-1}
}.
\end{align*}
Therefore, it follows that
\begin{align*}
P=L_1^{-1},\qquad
Q=-M^{-1}NL_1^{-1}.
\end{align*}
Because $L_1$ is positive definite, $P$ is positive definite.
Using a positive definite $Y$ such that $M^{-1}=YY$, we have
\[
Q= -Y 
\underbrace{
Y N Y
\left(
I + (Y N Y )^2
\right)^{-1}
}_{X}
Y.
\]
Because $N$ in \eqref{eq:defMNm} is negative semidefinite, $Y N Y$ can be diagonalized by an orthonormal matrix $V$.
Thus, we have
\begin{align}\label{eq:defXm}
X = V  \varLambda \left(
I + \varLambda^2
\right)^{-1}
V^{\sf T}
\end{align}
where $\varLambda$ is a negative semidefinite diagonal matrix composed of the eigenvalues of $Y N Y$.
Therefore, $Q$ is positive semidefinite because $X$ is negative semidefinite.
Finally, from 
\[
\bm{Y}^{\rm red} = \underbrace{Q}_{G^{\rm red}} + \bm{j} \underbrace{(-P)}_{B^{\rm red}},
\]
we see that $G^{\rm red}$ is positive semidefinite, and $B^{\rm red}$ is negative definite.
\end{IEEEproof}

Proposition~\ref{prop:defBGred} states that the definiteness of the conductance and susceptance matrices are invariant under the Kron reduction.
The following property of the reduced conductance matrix is also important.

\begin{proposition}\label{prop:GredG}
Consider the admittance matrix $\bm{Y}$ in \eqref{eq:GandB}, and the reduced admittance matrix $\bm{Y}^{\rm red}$ in \eqref{eq:Yred}.
Assume that the inequality condition in Lemma~\ref{lem:gamnon} is satisfied.
Then
\begin{equation}\label{eq:Gredker}
\sfspan \bigl\{ 
\sfdiag \left(
1- \beta_i X_i'
\right) \mathds{1} 
\bigr\}
\subseteq
\sfker G^{\rm red}.
\end{equation}
Furthermore, $G^{\rm red}$ is zero if and only if  $G$ is zero.
\end{proposition}

\begin{IEEEproof}
First, we prove \eqref{eq:Gredker}.
From \eqref{eq:kerGB0}, we see that
\[
\bm{\varGamma} 
\sfdiag
\left(
\tfrac{1}{X_i'}
\right)
\mathds{1}
=
\sfdiag
\bigl(
1 - \beta_i  X_i'
\bigr)
\mathds{1}.
\]
Multiplying it by $\bm{Y}^{\rm red}$ from the left side, we have
\[
-\bm{j} \sfdiag
\left(
\tfrac{1}{X_i'}
\right) \mathds{1}
=
\bigl(
G^{\rm red} + \bm{j} B^{\rm red} 
\bigr)
\sfdiag
\bigl(
1 - \beta_i  X_i'
\bigr)
\mathds{1}
\]
Considering the real part of both sides, we have
\[
0 = G^{\rm red} \sfdiag
\bigl(
1 - \beta_i  X_i'
\bigr)
\mathds{1},
\]
which proves \eqref{eq:Gredker}.

Next, we prove the second claim.
Consider the same notation as that in the proof of Proposition~\ref{prop:defBGred}.
Because $Y $ is positive definite, we see that
\[
G^{\rm red} =0 
\ \ 
\Longleftrightarrow
\ \ 
X =0
\ \ 
\Longleftrightarrow
\ \ 
N=0
\ \ 
\Longleftrightarrow
\ \ 
G=0
\]
where we have used the facts that $G^{\rm red}$ is equal to $Q$, $X$ is defined as in \eqref{eq:defXm}, and $N$ is defined as in \eqref{eq:defMNm}.
\end{IEEEproof}

Proposition~\ref{prop:GredG} shows the property of the reduced conductance matrix.
In particular, it is zero if and only if the original one before the Kron reduction is lossless.
%In particular, if the ground capacitances of all transmission lines are zero, then $\sfspan\{\mathds{1}\}$ is again involved  in its kernel.
A notable point in Propositions~\ref{prop:defBGred} and \ref{prop:GredG} is that the properties of the reduced admittance matrix are derived from those of the original admittance matrix before the Kron reduction.

\section{Positive Semidefiniteness of Synchronizing Torque Coefficient Matrix}\label{sec:EIclass}

We aim at deriving a sufficient condition for the Hessian $\nabla^2 \tilde{U}( \delta^{\star})$ in Theorem~\ref{thm:Hesconcl}, i.e., the synchronizing torque coefficient matrix, to be positive semidefinite.
Recall that this Hessian involves the elements of
\begin{align}\label{eq:Btilred}
\tilde{B}^{\rm red}= -
\left\{
\sfdiag ( X_i  )
- \sfdiag \left( X_i \right) B \sfdiag \left( X_i \right)
\right\}^{-1}.
\end{align}
As shown in the following lemma, all elements of this $\tilde{B}^{\rm red}$ are proven to be  non-positive.

\begin{lemma}\label{lem:Bsign}
Consider the susceptance matrix $B$ in \eqref{eq:defbeta}. 
For the constants $\beta_1,\ldots,\beta_N$ in \eqref{eq:defbeta}, if
\begin{align}\label{eq:Xsbeta}
\beta_i X_i \leq 1
,\quad \forall i \in \{1,\ldots,N\}
\end{align}
and the strict inequality holds for at least one bus, then $\tilde{B}^{\rm red}$ in \eqref{eq:Btilred} exists, and its all elements are non-positive.
\end{lemma}

\begin{IEEEproof}
Denote  $\tilde{B}^{\rm red}$ as $-K^{-1}$.
Then, we see that
\[
K= \sfdiag \bigl( X_i (1- \beta_i X_i ) \bigr) 
- \sfdiag \left( X_i \right) B_0 \sfdiag \left( X_i \right).
\]
Because $-B_0$ is a weighted graph Laplacian, $K$ is a positive definite matrix whose off-diagonal elements are non-positive if the condition in the claim is satisifed.
That is, $K$ is a nonsinglar M-matrix \cite{bernstein2009matrix}.
Therefore, all elements of $K^{-1}$ are non-negative.
This proves the claim.
\end{IEEEproof}

Based on Lemma~\ref{lem:Bsign}, we can find a simple sufficient condition for $\nabla^2 \tilde{U}( \delta^{\star})$ to be positive semidefinite as follows.

\begin{proposition}\label{prop:hessposd}
Consider the susceptance matrix $B$ in \eqref{eq:defbeta}, and the strain energy function $\tilde{U}$ in \eqref{eq:Udel}. 
Assume that the inequality condition in Lemma~\ref{lem:Bsign} is satisfied.
If 
\begin{equation}\label{eq:delrange}
\bigl| \delta_i^{\star} - \delta_j^{\star} \bigr| \leq \frac{\pi}{2}
,\quad
\forall (i,j)\in  \{1,\ldots,N\} \times \{1,\ldots,N\},
\end{equation}
then $\nabla^2 \tilde{U}( \delta^{\star})$ is positive semidefinite.
\end{proposition}

\begin{IEEEproof}
Because $\nabla^2 \tilde{U}( \delta^{\star})$ is a wighted graph Laplacian, it can be represented in the form
\[
\nabla^2 \tilde{U}( \delta^{\star})
=
-W
\sfdiag \Bigl( 
V_{{\rm fd}i}^{\star}V_{{\rm fd}j}^{\star} \tilde{B}^{\rm red}_{ij} \sfcos \delta_{ij}^{\star} 
\Bigr)_{(i,j) \in \mathcal{C}}
W^{\sf T}
\]
where $W \in \mathbb{R}^{N\times \frac{N(N-1)}{2}}$ denotes the incidence matrix associated with the complete graph, and
\[
\mathcal{C}:= \bigl\{
(i,j) \in  \{1,\ldots,N\} \times \{1,\ldots,N\} : i < j
\bigr\}.
\]
Note that every $\tilde{B}^{\rm red}_{ij}$ is less than or equal to 0 as shown in Lemma~\ref{lem:Bsign}.
Thus, the claim is proven.
\end{IEEEproof}

Proposition~\ref{prop:hessposd} states that if the differences among the steady-state rotor angles of all generators are small, then the classical generator network is EI passive, or equivalently, the synchronizing torque coefficient matrix is positive semidefinite.
This natural consequence is derived from the fact that the susceptance matrix before the Kron reduction is negative semidefinite.
\ifCLASSOPTIONcaptionsoff
  \newpage
\fi



%\bibliographystyle{IEEEtran}
%\begin{thebibliography}{1}


%\end{thebibliography}

\bibliographystyle{IEEEtran}% bib style
\bibliography{reference,reference_CREST}% your bib database

\begin{thebibliography}{10}
\providecommand{\url}[1]{#1}
\csname url@samestyle\endcsname
\providecommand{\newblock}{\relax}
\providecommand{\bibinfo}[2]{#2}
\providecommand{\BIBentrySTDinterwordspacing}{\spaceskip=0pt\relax}
\providecommand{\BIBentryALTinterwordstretchfactor}{4}
\providecommand{\BIBentryALTinterwordspacing}{\spaceskip=\fontdimen2\font plus
\BIBentryALTinterwordstretchfactor\fontdimen3\font minus
  \fontdimen4\font\relax}
\providecommand{\BIBforeignlanguage}[2]{{%
\expandafter\ifx\csname l@#1\endcsname\relax
\typeout{** WARNING: IEEEtran.bst: No hyphenation pattern has been}%
\typeout{** loaded for the language `#1'. Using the pattern for}%
\typeout{** the default language instead.}%
\else
\language=\csname l@#1\endcsname
\fi
#2}}
\providecommand{\BIBdecl}{\relax}
\BIBdecl

\bibitem{ishizaki2021necessity}
T.~Ishizaki and A.~Chakrabortty, ``Necessity of lossless transmission and
  convexity of potential energy function for equilibrium independent passivity
  of power systems,'' in \emph{Proc. of The 60th IEEE Conference on Decision
  and Control}, 2021.

\bibitem{kundur1994power}
P.~Kundur, \emph{Power system stability and control}.\hskip 1em plus 0.5em
  minus 0.4em\relax Tata McGraw-Hill Education, 1994.

\bibitem{arcak2007passivity}
M.~Arcak, ``Passivity as a design tool for group coordination,''
  \emph{Automatic Control, IEEE Transactions on}, vol.~52, no.~8, pp.
  1380--1390, 2007.

\bibitem{yang2019distributed}
P.~Yang, F.~Liu, Z.~Wang, and C.~Shen, ``Distributed stability conditions for
  power systems with heterogeneous nonlinear bus dynamics,'' \emph{IEEE
  Transactions on Power Systems}, 2019.

\bibitem{trip2016internal}
S.~Trip, M.~B{\"u}rger, and C.~De~Persis, ``An internal model approach to
  (optimal) frequency regulation in power grids with time-varying voltages,''
  \emph{Automatica}, vol.~64, pp. 240--253, 2016.

\bibitem{stegink2017unifying}
T.~Stegink, C.~De~Persis, and A.~van~der Schaft, ``A unifying energy-based
  approach to stability of power grids with market dynamics,'' \emph{Automatic
  Control, IEEE Transactions on}, vol.~62, no.~6, pp. 2612--2622, 2017.

\bibitem{wang2018distributed}
Z.~Wang, F.~Liu, J.~Z. Pang, S.~H. Low, and S.~Mei, ``Distributed optimal
  frequency control considering a nonlinear network-preserving model,''
  \emph{IEEE Transactions on Power Systems}, vol.~34, no.~1, pp. 76--86, 2018.

\bibitem{li2016connecting}
N.~Li, L.~Chen, C.~Zhao, and S.~H. Low, ``Connecting automatic generation
  control and economic dispatch from an optimization view,'' \emph{Control of
  Network Systems, IEEE Transactions on}, vol.~3, no.~3, pp. 254--264, 2015.

\bibitem{ortega2005transient}
R.~Ortega, M.~Galaz, A.~Astolfi, Y.~Sun, and T.~Shen, ``Transient stabilization
  of multimachine power systems with nontrivial transfer conductances,''
  \emph{Automatic Control, IEEE Transactions on}, vol.~50, no.~1, pp. 60--75,
  2005.

\bibitem{sauer2017power}
P.~W. Sauer, M.~A. Pai, and J.~H. Chow, \emph{Power system dynamics and
  stability: with synchrophasor measurement and power system toolbox}.\hskip
  1em plus 0.5em minus 0.4em\relax John Wiley \& Sons, 2017.

\bibitem{sadamoto2019dynamic}
T.~Sadamoto, A.~Chakrabortty, T.~Ishizaki, and J.-i. Imura, ``Dynamic modeling,
  stability, and control of power systems with distributed energy resources:
  Handling faults using two control methods in tandem,'' \emph{IEEE Control
  Systems Magazine}, vol.~39, no.~2, pp. 34--65, 2019.

\bibitem{hines2011equilibrium}
G.~H. Hines, M.~Arcak, and A.~K. Packard, ``Equilibrium-independent passivity:
  A new definition and numerical certification,'' \emph{Automatica}, vol.~47,
  no.~9, pp. 1949--1956, 2011.

\bibitem{simpson2019equilibrium}
J.~W. Simpson-Porco, ``Equilibrium-independent dissipativity with quadratic
  supply rates,'' \emph{IEEE Transactions on Automatic Control}, vol.~64,
  no.~4, pp. 1440--1455, 2019.

\bibitem{DePersis2017}
C.~De~Persis and N.~Monshizadeh, ``Bregman storage functions for microgrid
  control,'' \emph{IEEE Transactions on Automatic Control}, vol.~63, no.~1, pp.
  53--68, 2017.

\bibitem{chang1995direct}
H.-D. Chang, C.-C. Chu, and G.~Cauley, ``Direct stability analysis of electric
  power systems using energy functions: theory, applications, and
  perspective,'' \emph{Proceedings of the IEEE}, vol.~83, no.~11, pp.
  1497--1529, 1995.

\bibitem{narasimhamurthi1984existence}
N.~Narasimhamurthi, ``On the existence of energy function for power systems
  with transmission losses,'' \emph{IEEE transactions on circuits and systems},
  vol.~31, no.~2, pp. 199--203, 1984.

\bibitem{Anderson2008}
P.~M. Anderson and A.~A. Fouad, \emph{Power system control and
  stability}.\hskip 1em plus 0.5em minus 0.4em\relax John Wiley \& Sons, 2008.

\bibitem{ishizaki2018graph}
T.~Ishizaki, A.~Chakrabortty, and J.-i. Imura, ``Graph-theoretic analysis of
  power systems,'' \emph{Proceedings of the IEEE}, vol. 106, no.~5, pp.
  931--952, 2018.

\bibitem{simpson2013synchronization}
J.~W. Simpson-Porco, F.~D{\"o}rfler, and F.~Bullo, ``Synchronization and power
  sharing for droop-controlled inverters in islanded microgrids,''
  \emph{Automatica}, vol.~49, no.~9, pp. 2603--2611, 2013.

\bibitem{guggilam2017optimizing}
S.~S. Guggilam, C.~Zhao, E.~Dall’Anese, Y.~C. Chen, and S.~V. Dhople,
  ``Optimizing power--frequency droop characteristics of distributed energy
  resources,'' \emph{IEEE Transactions on Power Systems}, vol.~33, no.~3, pp.
  3076--3086, 2017.

\bibitem{monshizadeh2019conditions}
N.~Monshizadeh, P.~Monshizadeh, R.~Ortega, and A.~van~der Schaft, ``Conditions
  on shifted passivity of port-hamiltonian systems,'' \emph{Systems \& Control
  Letters}, vol. 123, pp. 55--61, 2019.

\bibitem{rockafellar1970convex}
R.~T. Rockafellar, \emph{Convex analysis}.\hskip 1em plus 0.5em minus
  0.4em\relax Princeton university press, 1970.

\bibitem{petersen2010feedback}
I.~R. Petersen and A.~Lanzon, ``Feedback control of negative-imaginary
  systems,'' \emph{IEEE Control Systems Magazine}, vol.~30, no.~5, pp. 54--72,
  2010.

\bibitem{chen2020modelling}
M.~Chen, D.~Zhou, and F.~Blaabjerg, ``Modelling, implementation, and assessment
  of virtual synchronous generator in power systems,'' \emph{Journal of Modern
  Power Systems and Clean Energy}, vol.~8, no.~3, pp. 399--411, 2020.

\bibitem{kron1939tensor}
G.~Kron, \emph{Tensor Analysis of Networks}.\hskip 1em plus 0.5em minus
  0.4em\relax New York: Wiley, 1939.

\bibitem{bernstein2009matrix}
D.~S. Bernstein, \emph{Matrix mathematics: theory, facts, and formulas}.\hskip
  1em plus 0.5em minus 0.4em\relax Princeton University Press, 2009.

\end{thebibliography}
%\bibliography{sshrrefs}% your bib database

\end{document}